\pdfoutput=1
\documentclass[preprint,11pt]{elsarticle}
\usepackage{array}
\usepackage{arydshln}
\usepackage{amsmath,amsfonts,booktabs,amssymb}
\usepackage{tikz,pgfplots}
\usepackage{graphicx,subfigure}
\usepackage{placeins}
\usepackage{multirow}
\usepackage{enumerate}
\usepackage{todonotes}
\usepackage{hyperref}
\usepackage{textcomp,gensymb,url}
\usepackage[left=1in, right=1.5in, top=1in, bottom=1in]{geometry}
\usetikzlibrary{external} 

\journal{Journal of Computational Physics}
\pgfplotsset{mystyle/.style={mark=none}}
\begin{document}
\begin{frontmatter}
\title{A Fast Block Low-Rank Dense Solver with Applications to Finite-Element Matrices} 


\author[Amir]{AmirHossein Aminfar\corref{cor1}\fnref{label1}}
\ead{aminfar@stanford.edu}
\author[Siva]{Sivaram Ambikasaran\fnref{label2}} 
\author[Eric]{Eric Darve\fnref{label1}}

\fntext[label1]{Mechanical Engineering Department, Stanford University}
\cortext[cor1]{Corresponding author. +1 650-644-7624}
\fntext[label2]{Courant Institute of Mathematical Sciences, New York University}
\address[Amir]{496 Lomita Mall, Room 104, Stanford, CA, 94305}
\address[Siva]{Warren Weaver Hall, Room-1105A, 251, Mercer Street, New York, NY 10012}
\address[Eric]{496 Lomita Mall, Room 209, Stanford, CA, 94305}

\begin{abstract}

This article presents a fast solver for the dense ``frontal'' matrices that arise from the multifrontal sparse elimination process of $3$D elliptic PDEs. The solver relies on the fact that these matrices can be efficiently represented as a hierarchically off-diagonal low-rank (HODLR) matrix. To construct the low-rank approximation of the off-diagonal blocks, we propose a new pseudo-skeleton scheme, the boundary distance low-rank approximation, that picks rows and columns based on the location of their corresponding vertices in the sparse matrix graph. We compare this new low-rank approximation method to the adaptive cross approximation (ACA) algorithm and show that it achieves betters speedup specially for unstructured meshes. Using the HODLR direct solver as a preconditioner (with a low tolerance) to the GMRES iterative scheme, we can reach machine accuracy much faster than a conventional LU solver. Numerical benchmarks are provided for frontal matrices arising from $3$D finite element problems corresponding to a wide range of applications.
\end{abstract}


\end{frontmatter}
\section{Introduction}
\label {sec:intro}
In many engineering applications, solving large finite element systems is of great significance.
Consider the system 
\[
Ax = b
\]
arising from the finite element discretization of an elliptic PDE, where $A\in \mathbb{R}^{N \times N}$ is a sparse matrix with a symmetric pattern. In many practical applications, the matrix $A$ might be ill-conditioned and thus, challenging for iterative methods. On the other hand, conventional direct solver algorithms, while being robust in handling ill-conditioned matrices, are computationally expensive ($\mathcal{O}(N^{1.5})$ for 2D meshes and $\mathcal{O}(N^{2})$ for 3D meshes). Since one of the main bottlenecks in the direct multifrontal solve process is the high computational cost of solving dense frontal matrices, we mainly focus on solving these matrices in this article. Our goal is to build an iterative solver, which utilizes a fast direct solver as a preconditioner for the dense frontal matrices. The direct solver in this scheme acts as a highly accurate pre-conditioner. This approach combines the advantages of the iterative and direct solve algorithms, i.e., it is fast, accurate and robust in handling ill-conditioned matrices.

To be consistent with our previous work, we adopt the notations used in~\cite{SivaFDS}. We should also mention that `$n$' refers to the size of dense matrices and `$N$' refers to the size of sparse matrices (e.g., number of degrees of freedom in a finite-element mesh).

In the next section, we review the previous literature on both dense structured solvers and sparse multifrontal solvers. We then introduce a hierarchical off-diagonal low-rank (from now on abbreviated as HODLR) direct solver in Section~\ref{sec:directSolver}. In Section~\ref{sec:introLR}, we introduce the boundary distance low-rank (BDLR) algorithm as a robust low-rank approximation scheme for representing the off-diagonal blocks of the frontal matrices. Section~\ref{sec:MF} discusses the application of the iterative solver with a fast HODLR direct solver preconditioner to the sparse multifrontal solve process and demonstrates the solver for a variety of $3$D meshes. We also show an application in combination with the FETI-DP method~\cite{FETI_DP1}, which is a family of domain decomposition algorithms to accelerate finite-element analysis on parallel computers. We present the results and numerical benchmarks in Section~\ref{sec:benchmarks}.

\section{Previous Work}
\label{sec:prevWork}

\subsection{Fast Direct Solvers for Dense Hierarchical Matrices}
\label{sec:denseMatrices}
Hierarchical matrices are data sparse representation of a certain class of dense matrices. This representation relies on the fact that these matrices can be sub-divided into a hierarchy of smaller block matrices, and certain sub-blocks (based on the admissibility criterion) can be efficiently represented as a low-rank matrix. We refer the readers to~\cite{hackbusch1999sparse, hackbusch2000sparse, grasedyck2003construction, hackbusch2002data, borm2003hierarchical,  ULV, chandrasekaran2006fast1} for more details. These matrices were introduced in the context of integral equations~\cite{hackbusch1999sparse,hackbusch2000sparse,ying2009fast,lai2014fast} arising out of elliptic partial differential equations and potential theory. Subsequently, it has also been observed that dense fill-ins in finite element matrices~\cite{xia2009superfast}, radial basis function interpolation~\cite{SivaFDS}, kernel density estimation in machine learning, covariance structure in statistic models~\cite{chen2014data}, Bayesian inversion~\cite{SivaFDS, ambikasaran2013large, ambikasaran2013fastBayes}, Kalman filtering~\cite{li2014kalman}, and Gaussian processes~\cite{ambikasaran2014fastdet}, can also be efficiently represented as data-sparse hierarchical matrices. Broadly speaking, these matrices can be grouped into two general categories based on the admissibility criterion: (i) Strong admissibility: sub-block that correspond to the interaction between well-separated clusters are low-rank; (ii) Weak admissibility: sub-block corresponding to non-overlapping interactions are low-rank. Ambikasaran~\cite{ambikasaran2013thesis} provides a detailed description of these different hierarchical structures.

We review some of the previously developed structured dense solvers for hierarchical matrices and discuss them in relation to our work. Hackbusch~\cite{hackbusch1999sparse,grasedyck2003construction} introduced the concept of $\mathcal{H}$-matrices, which are the most general class of hierarchical matrices with the strong admissibility criterion~\cite{hackbusch1999sparse,grasedyck2003construction,hackbusch2002data,hackbusch2000h2,hackbusch2002h2,hackbusch2000sparse,hackbusch2000sparse1,bebendorf2000approximation,bebendorf2008hierarchical,borm2003hierarchical}. Contrary to the HODLR matrix structure, in which the off-diagonal blocks are low-rank, in $\mathcal{H}$-matrices, the off-diagonal blocks are further decomposed into low-rank and full-rank blocks. Thus, the rank can be kept small. In HODLR, we make a single low-rank approximation for the off-diagonal blocks and the rank is larger as a result. Hence, the HODLR structure makes for a much simpler representation and is often used because of its simplicity compared to the $\mathcal{H}$-matrix structure. Hackbusch~\cite{grasedyck2003construction} suggests a recursive block low-rank factorization scheme for $\mathcal{H}$-matrices. This method is based on the idea that all the dense matrix algebra (matrix multiplication and matrix addition) can be replaced by $\mathcal{H}$-matrix algebra. As a result, the inverse of an $\mathcal{H}$-matrix can also be approximated as an $\mathcal{H}$-matrix itself. This results in a computational complexity of  $\mathcal{O}(n\log^2(n))$ for an $\mathcal{H}$-matrix factorization.

We note that the approach in this paper is based on the Woodbury matrix identity. It is therefore different from the algorithm in Hackbusch~\cite{grasedyck2003construction} for example. The latter is based on a block LU factorization, while the Woodbury identity reduces the global solve to block diagonal solves followed by a correction update.

The HODLR matrix structure is the most general off-diagonal low-rank structure with weak admissibility. Solvers for this matrix class have a computational cost of $\mathcal{O}(n\log^2n)$. In an HODLR matrix, the off-diagonal low-rank bases do not have a nested structure across different levels~\cite{SivaFDS}. The HSS matrix is an HODLR matrix but, in addition, has a nested off-diagonal low-rank structure. Solvers for the HSS matrices have an $\mathcal{O}(n)$ complexity~\cite{FastHSS,HSSSuperFastMF}.

Martinsson and Rokhlin~\cite{MartinssonRokhlin} discuss an $\mathcal{O}(n)$ direct solver for boundary integral equations based on the HSS structure. Their method is based on the fact that for a matrix of rank $r$, there exists a well-conditioned column operation, which leaves $r$ columns unchanged and sets the remaining columns to zero. Using this idea, they derive a two-sided compressed factorization of the inverse of the HSS matrix. Their generic algorithm requires $\mathcal{O}(n^2)$ operations to construct the inverse. However, they accelerate their algorithm to $\mathcal{O}(n \log^{\kappa}(n))$ when applied to two-dimensional contour integral equations.

Chandrasekaran et al.~\cite{ULV} present a fast $\mathcal{O}(n)$ direct solver for HSS matrices. In their article, they construct an implicit $ULV^H$ factorization of an HSS matrix, where $U$ and $V$ are unitary matrices, $L$ is a lower triangular matrix and $^H$ is the transpose conjugate operator. Their method is based on the Woodbury matrix identity and the fact that for a low-rank representation of the form $UBV^H$, where $U$ and $V$ are thin matrices with $r$ columns, there exists a unitary transformation $Q$, in which only the last $r$ rows of $QU$ are nonzero. They use this observation to recursively solve the linear system of equations. Since this method requires constructing an HSS tree, the authors suggest an algorithm that uses the SVD or the rank revealing QR decomposition, recursively, to construct the HSS tree in $\mathcal{O}(n^2)$ time.

Gillman et al.~\cite{GillmanFDS} discuss an $\mathcal{O}(n)$ algorithm for directly solving integral equations in one-dimensional domains. The algorithm relies on applying the Sherman-Morrison-Woodbury formula (see for example~\cite{SivaFDS}) recursively to an HSS tree structure to achieve $\mathcal{O}(r^2n)$ solve complexity, where $r$ is the rank of the off diagonal blocks in the HSS matrix. They also describe an $\mathcal{O}(r^2n)$ algorithm for constructing an HSS representation of the matrix resulting from a Nystr\"om discretization of a boundary integral equation.

Ho and Greengard~\cite{HoFDS} present a fast direct solver for HSS matrices. They use the interpolative decomposition (ID) (see for example~\cite{IDref}) algorithm with random sampling to obtain the low-rank representations of the off-diagonal blocks. The computational complexity of the low-rank approximation algorithm is $\mathcal{O}(mn\log r + r^2n)$ for a matrix $K\in \mathbb{R}^{m \times n}$. After obtaining the hierarchical matrix representation of the original dense matrix, new variables and equations are introduced into the system of equations. Finally, all equations are assembled into an extended sparse matrix and a conventional sparse solver is used to factorize the sparse matrix. This method has a computational complexity of $\mathcal{O}(n)$ for both the pre-computation and solution phases for boundary integral equations in 2D, while in 3D, these phases cost $\mathcal{O}(n^{1.5})$ and $\mathcal{O}(n\log(n))$ respectively.

Kong et al.~\cite{KongFDS} have developed an $\mathcal{O}(n^2)$ dense solver for HODLR matrices. Similar to~\cite{MartinssonRokhlin}, they accelerate their algorithm to $\mathcal{O}(n\log^2(n))$, when applied to boundary integral equations. Their method uses the Sherman-Morrison-Woodbury formula to construct a one-sided hierarchical factorization of the inverse of these matrices, in which each factor is a block diagonal matrix. The low-rank approximation scheme in their paper is based on the rank revealing QR algorithm. The authors use the pivoted Gram-Schmidt algorithm to obtain $r$ orthogonal basis vectors for the low-rank matrix in question. For a matrix $K\in \mathbb{R}^{m\times n}$ with rank $r$, this low-rank approximation scheme requires $\mathcal{O}(mnr)$ operations. Then, they use a randomized algorithm from~\cite{randomized} to accelerate their low-rank approximation scheme. This accelerated low-rank approximation algorithm costs $\mathcal{O}(mn \, \log (l+lnr))$ in the general case where $r<l<\min(m,n)$. 

Ambikasaran and Darve~\cite{SivaFDS} present an $\mathcal{O}(n\log^2(n))$ solver for HODLR matrices and an $\mathcal{O}(n\log(n))$ solver for p-HSS matrices. This approach differs from the approach mentioned in~\cite{KongFDS} in the fact that, while~\cite{KongFDS} constructs the inverse,~\cite{SivaFDS} constructs a factorization of the matrix. Each factor in this factorization scheme is a block diagonal matrix with each block being a low-rank perturbation of the identity matrix. The authors then use the Sherman-Morrison-Woodbury formula to invert each block in the block diagonal factors. The article uses the Chebyshev low-rank approximation scheme to factorize the off-diagonal blocks.

As mentioned above, solvers for the HSS matrix structure have the lowest computational complexity \textemdash{} $\mathcal{O}(r^2n)$, $r$ being the rank of approximation \textemdash{} among other hierarchically off-diagonal low-rank matrix structures. While the HSS structure is attractive, the nested structure makes it more complicated and more difficult to work with, compared to the simpler HODLR structure. Furthermore, the off-diagonal rank increases from root to leaves in the HSS tree, whereas the off-diagonal ranks at each level are independent from each other in the HODLR structure. This often leads to lower average off-diagonal rank in the HODLR structure compared to HSS.



A point worth mentioning is that the solver discussed in the current article relies on purely algebraic technique (instead of analytic or geometry based techniques) to construct the low-rank approximation of the off-diagonal blocks. Analytic low-rank approximation techniques like the Chebyshev low-rank approximation, multipole expansions, etc., are only applicable when the matrix definition involves an analytical kernel function. 
In this article, we propose a boundary distance low-rank approximation (from now on abbreviated as BDLR), which relies on the underlying sparse matrix graph to choose the desired rows and columns in constructing a low-rank representation. We also compare with the adaptive cross approximation algorithm~\cite{ACA} (from now on abbreviated as ACA), which is also a purely algebraic scheme to construct low-rank approximations of the off-diagonal blocks. 

Due to its black-box nature, the solver can handle a wide range of dense matrices arising from boundary integral equations, covariance matrices in statistics, frontal matrices arising in the context of finite-element matrices, etc. Table~\ref{table:summaryHODLR} summarizes the dense solver algorithms mentioned above.

\begin{table}[htbp]
\centering
\scalebox{0.9}{
\begin{tabular}{p{2.8cm}p{2.3cm}p{4.6cm}p{5cm}}
	\toprule
	\small{Article} & \small{Matrix Class} & \small{Factorization} & \small{Application}\\ \midrule
	\small{Hackbusch}~\cite{hackbusch1999sparse,grasedyck2003construction}&\small{ $\mathcal{H}$}& Recursive block factorization of the matrix& \small{BEM integral operators}\\ \midrule
	\small{Martinsson and Rokhlin~\cite{MartinssonRokhlin}} & \small{HSS} & \small{Two sided compressed factorization of the inverse} & \small{2D boundary integral equations}\\ \midrule
	\small{Chandrasekaran et al.~\cite{ULV}} & \small{HSS} & \small{$ULV^H$ factorization of the matrix} & \small{Radial basis function matrices}\\ \midrule
	\small{Gillman et al.~\cite{GillmanFDS}} & \small{HSS} & \small{Data sparse factorization of the inverse} & \small{1D integral equations with Nystr\"om discretization}\\ \midrule
	\small{Ho and Greengard~\cite{HoFDS}} & \small{HSS} & \small{Factorization of the extended sparse system} & \small{2D and 3D boundary integral equations}\\ \midrule
	\small{Kong et al.~\cite{KongFDS}} & \small{HODLR} & \small{One sided hierarchical factorization of the inverse} & \small{Boundary integral equations}\\ \midrule
	\small{Ambikasaran and Darve~\cite{SivaFDS}} & \small{HODLR, p-HSS} & \small{Block-diagonal factorization of the matrix} & \small{Interpolation using radial basis functions}\\ \midrule
	\small{This article} & \small{HODLR} & \small{Recursive block LU factorization of the matrix}& \small{Finite-element matrices}\\ \bottomrule
\end{tabular}
}
\caption{Summary of fast dense structured solvers.}
\label{table:summaryHODLR}
\end{table}

\subsection{Fast Direct Solvers for Sparse Matrices}
\label{sec:sparseMatrices}
As mentioned in Section~\ref{sec:denseMatrices}, we are interested in accelerating the direct solve process for finite-element matrices. In this article, we focus on the finite-element discretization of elliptic PDEs. One common way of factorizing such matrices is using the sparse Cholesky factorization. The efficiency of this algorithm strongly depends on the ordering of mesh nodes~\cite{MFGeneralMesh}. Sparse Cholesky factorization takes $\mathcal{O}(N^2)$ flops in 2D with a typical row-wise or column-wise mesh ordering, where $N$ is the number of degrees of freedom~\cite{xia2009superfast}. The most efficient method for solving such matrices is the multifrontal method with nested dissection~\cite{NestedDissection}, which takes $\mathcal{O}(N^{1.5})$ flops for two-dimensional and $\mathcal{O}(N^{2})$ for three dimensional meshes~\cite{MFGeneralMesh}.

The multifrontal method was originally introduced by Duff \& Reid~\cite{originalMF}, George~\cite{NestedDissection} and Liu~\cite{LiuMF}, as an extension to the frontal method of Irons \cite{originalFrontal}. In this algorithm, the overall factorization is done by factorizing smaller dense frontal matrices~\cite {MFReview}. For each node or super-node in the elimination tree, the frontal matrix is obtained using an update process called the ``extend-add" process, which involves updates from the previously eliminated nodes.

Martinsson~\cite{SpiralElimination} uses a spiral elimination approach along with HSS compression of Schur complements to achieve $\mathcal{O}(N\log^2N)$ time complexity. This approach is not based on the multifrontal method and requires a mesh that can be partitioned into concentric annuli.

Gillman et al.~\cite{HBSGillman} proposed an accelerated nested dissection algorithm for obtaining the Dirichlet-to-Neumann operator associated with a 2D elliptic boundary value problem. The authors approximate the Schur complements that appear in the elimination process as hierarchically block separable (HBS) matrices, a structure similar to HSS matrices. Using this matrix structure, they are able to obtain the Dirichlet-to-Neumann operator with a cost of $\mathcal{O}(N)$ compared to $\mathcal{O}(N^{1.5})$ of the conventional multifrontal method with nested dissection.

There have been some recent efforts to reduce the computational cost of the multifrontal method with nested dissection. Xia et al.~\cite{xia2009superfast} observed that frontal and update matrices in the multifrontal elimination process can be approximated with hierarchically semi-separable (HSS) matrices. The authors develop a structured extend-add process to facilitate the formation of the frontal matrices using the HSS structure. Next, they perform a structured dense Cholesky factorization on the obtained frontal matrix. The authors use the algorithm in~\cite{HSSSuperFastMF} to compute the explicit factorizations of HSS matrices. Using this procedure, they are able to achieve nearly linear time complexity for $2$D meshes. However, only regular well shaped meshes in $2$D are considered in the article. Schmitz et al.~\cite {MFGeneralMesh} extend the approach of~\cite{xia2009superfast} to a more general setting of unstructured and adaptive grids in 2D. 

Xia~\cite{GeneralizedMF} introduced an efficient multifrontal factorization for general sparse matrices. The author approximates the frontal matrices using the HSS structure and introduces the concept of reduced HSS matrices that reduce the computational cost of operation on HSS matrices. For simplicity, this approach keeps the update matrices as dense matrices which leads to high memory consumption for large sparse matrices. Xia~\cite{randomizedMF} introduces a new algorithm that overcomes this deficiency by randomization. That is, instead of passing dense update matrices along the elimination tree, this approach passes a skinny randomized matrix vector product. In addition to saving memory, this approach only requires skinny extend adds (extend adds on all rows and only a subset of columns) which leads to improvements in efficiency. This method is based on the work of Martinsson~\cite{RandomizedHSS} which provides an algorithm for constructing HSS matrices using randomized matrix vector products.

Amestoy et al.~\cite{BLRMF} introduce a new low-rank matrix format called the Block Low-Rank (BLR) structure, a flat, non-hierarchical block matrix structure, for representing frontal matrices obtained in the multifrontal elimination process. The authors show that BLR is a good alternative to hierarchical structures like $\mathcal{H}$ and HSS matrices in terms of storage costs, flop count and parallelization for representing frontal matrices. The article demonstrates that the BLR format reduces the flop count and storage requirements for factorizing frontal matrices arising from a variety of large matrices coming from different physics applications. However, there is no discussion of the extend-add operations for BLR matrices. Furthermore, the article does not demonstrate a full multifrontal solver based on the BLR frontal matrix representation.

The approach presented here is based on the multifrontal method~\cite{MFReview}. It does not require constructing and maintaining HSS trees and can be applied to any mesh structure. Our method is based on the observation that the frontal matrices obtained during the multifrontal elimination process have an HODLR structure. This observation was also made by~\cite{xia2009superfast}. 

In order to factorize (eliminate) these frontal matrices, we present a dense HODLR structured solver. If the rank $r$ is $\mathcal{O}(1)$ (that is function of $\epsilon$ only), the algorithm has a computational cost of $\mathcal{O}(r^2 n\log^2 n)$ for an $n\times n$ frontal matrix. When solving 3D PDEs, we typically have that $r \in \mathcal{O}(n^{1/2})$. In that case, the computational cost is $\mathcal{O}(r^2 n)$, where $r$ is the largest rank found, at the top of the tree. This cost is, in fact, slightly favorable compared to what is reported for HSS in~\cite{randomizedMF} (see Table 4.3, p.\ 219), at least asymptotically for $n \to \infty$. The $\log^2 n$ factor disappears because the rank is bounded by a geometric series associated with the rank. 

We will benchmark the structured elimination (solve) process for frontal matrices corresponding to separators at various levels of the sparse elimination tree, for many different types of sparse matrices. It is worth mentioning that contrary to previous works which have mainly benchmarked matrices in the University of Florida Sparse Matrix Collections~\cite{UFLSparse}, we focus on frontal matrices arising from large and complicated mesh structures. These matrices are often very ill-conditioned and cannot be solved using traditional iterative techniques like GMRES~\cite{saad1986gmres} with diagonal preconditioning.  Our benchmarks show that obtaining a good preconditioner for unstructured meshes is significantly harder compared to structured meshes. Furthermore, solving 3D problems is an order of magnitude more difficult than 2D problems as the off-diagonal rank is significantly higher in 3D. Hence, this article mainly focuses on 3D meshes.

Table~\ref{table:summaryMF} shows a summary of various fast sparse matrix solvers in the literature.

\begin{table}[htbp]
\centering
\scalebox{0.8}{
\begin{tabular}{p{3.2cm}p{7cm}p{6cm}}
	\toprule
	\small{Article} & \small{Methodology} & \small{Test Cases \& Application} \\ \midrule
	Martinsson~\cite{SpiralElimination} & HSS compression and spiral elimination & Meshes that can be partitioned into concentric annuli \\ \midrule
	Gillman et al.~\cite{HBSGillman}&Approximating Schur complements as HBS matrices and using HBS algebra. & 2D elliptic boundary value problems discretized using a 5 point stencil on a regular square grid.  \\ \midrule
	Xia et al.~\cite{xia2009superfast} & HSS approximation of frontal matrices and structured extend-add& 2D structured meshes \\ \midrule
	Schmitz et al.~\cite{MFGeneralMesh} & Modified ~\cite{xia2009superfast} to accommodate adaptive and unstructured grids & 2D adaptive and unstructured meshes that roughly follow the pattern of a regular mesh \\ \midrule
	Xia~\cite{GeneralizedMF} & Introduction of reduced HSS matrices that reduce the operation cost on HSS matrices. For simplifications, the update matrices are kept as dense matrices. & Helmholtz Equation in 2D and University of Florida Sparse Matrix Collections~\cite{UFLSparse}\\ \midrule
	Xia~\cite{randomizedMF} & HSS compression using randomization techniques in~\cite{RandomizedHSS}. Passing randomized matrix vector products instead of dense update matrices and performing skinny extend-add operations. & Helmholtz Equation in 2D and University of Florida Sparse Matrix Collections~\cite{UFLSparse}\\ \midrule
	Amestoy et al.~\cite{BLRMF} & BLR format for representing frontal matrices. No discussion of BLR extend-add process. & Large matrices coming from different physics applications \\ \midrule
	
	\end{tabular}
}
\caption{Summary of fast sparse direct solvers.}
\label{table:summaryMF}
\end{table}

\section {An Iterative Solver with Direct Solver Preconditioning}
\label{sec:directIterative}
In this paper, we investigate using a fast HODLR direct solver as a preconditioner to the GMRES~\cite{saad1986gmres} iterative scheme. In this case, we use a relatively low accuracy for the direct solver. 

We will show that this approach is much faster than both a conventional LU solver and a high accuracy direct HODLR solver. We should also mention that this preconditioning method can be applied to any iterative solver (conjugate gradient (CG)~\cite{hestenes1952methods}, etc..).

\section {A Fast Direct Solver for HODLR Matrices}
One bottleneck of sparse solvers is the factorization of the dense frontal matrices that appear during the multifrontal elimination process. To accelerate the factorization of dense frontal matrices, we approximate them as HODLR matrices. As mentioned in Section~\ref{sec:denseMatrices}, HODLR matrices can be factorized in $\mathcal{O}(n\log^2n)$ which is a significant improvement over conventional dense factorizations which typically scale as $\mathcal{O}(n^3)$.

\label{sec:directSolver}
\subsection{HODLR Matrices}

A HODLR matrix has low-rank off-diagonal blocks at multiple levels.
As described in \cite{SivaFDS}, a 2-level HODLR matrix, $K\in \mathbb{R}^{n\times n}$, can be written as shown in Equation~\eqref{eq:HODLR2}:
\begin{align}
	K & =
	\begin{bmatrix}
		K_1^{(1)}&U_1^{(1)} (V_{1,2}^{(1)})^T \\
		U_2^{(1)} (V_{2,1}^{(1)})^T&K_2^{(1)}
	\end{bmatrix} \\
	& =
	\begin{bmatrix}
		\begin{bmatrix}
			K_1^{(2)}&U_1^{(2)} (V_{1,2}^{(2)})^T \\
			U_2^{(2)} (V_{2,1}^{(2)})^T&K_2^{(2)}
		\end{bmatrix}&
		         U_1^{(1)} (V_{1,2}^{(1)})^T \\
			U_2^{(1)} (V_{2,1}^{(1)})^T&
		\begin{bmatrix}
			K_3^{(2)}& (U_3^{(2)})^T (V_{3,4}^{(2)})^T \\
			U_4^{(2)} (V_{4,3}^{(1)})^T&K_4^{(2)}
		\end{bmatrix}
	\end{bmatrix}
\label{eq:HODLR2}	
\end{align}
where for a $p$-level HODLR matrix, $K_i^{(p)} \in \mathbb{R}^{n/2^p\times n/2^p}$, $U_{2i-1}^{(p)} , U_{2i}^{(p)}, V_{2i-1,2i}^{(p)}, V_{2i,2i-1}^{(p)} \in \mathbb{R}^{n/2^p\times r}$ and $r\ll n$. Further nested compression of the off-diagonal blocks will lead to HSS structures~\cite{SivaFDS}.

\subsection{Solver Derivation and Algorithm}
Contrary to the method introduced by Hackbusch~\cite{hackbusch1999sparse} which utilizes sequential block LU factorization, the HODLR direct solve algorithm presented in this section is based on the Woodbury matrix identity (see for example~\cite{hager1989updating,SivaFDS}). Although we do not use the formula explicitly, we perform the exact same operations. Looking at Equation~\eqref{eq:eqHODLR}, our method assumes that both diagonal blocks are nonsingular and factorizes them independently. However, Hackbusch~\cite{hackbusch1999sparse} only assumes that top diagonal block is invertible and factorizes the top diagonal block first. He then constructs the remaining Schur complement and continues on with the factorization. In comparing the two methods, one can see that because of the independent factorization of the diagonal blocks, the method presented in this section is better suited to parallel implementations.

Consider the following linear equation:
\begin{equation}
	\label{eq:denseSystem}
	Kx=F
\end{equation}	
where $K\in \mathbb{R}^{n\times n}$ is an HODLR matrix and $x$, $F\in \mathbb{R}^{n\times s}$. Now let's write $K$ as a one-level HODLR matrix and rewrite Equation~\eqref{eq:denseSystem} :
\begin{equation}
\label{eq:eqHODLR}
	K=
	\begin{bmatrix}
		K_1^{(1)}&U_1^{(1)} V_{1,2}^{(1)}{}^T \\
		U_2^{(1)} V_{2,1}^{(1)}{}^T&K_2^{(1)}
	\end{bmatrix}
	\begin{bmatrix}
		x_1^{(1)}\\
		x_2^{(1)}\\
	\end{bmatrix}
	=
	\begin{bmatrix}
		F_1\\
		F_2\\
	\end{bmatrix}
\end{equation}		
where $x_i^{(1)}, F_i^{(1)} \in \mathbb{R}^{(\frac{n}{2}\times s)}$. We now introduce two new variables $y_1^{(1)}$ and $y_2^{(1)}$:
\begin{align} 
	y_1^{(1)}&= V_{2,1}^{(1)}{}^Tx_1^{(1)} \\
	y_2^{(1)}&= V_{1,2}^{(1)}{}^Tx_2^{(1)} 
\end{align}
Rearranging~\eqref{eq:eqHODLR}, we have:
\begin{equation}
\label{eq:modifiedSystem}
	\underbrace{
		\begin{bmatrix}
			K_1^{(1)}&0&0&U_1^{(1)}\\
			0&K_2^{(1)}&U_2^{(1)}&0\\
			-V_{2,1}^{(1)}{}^T&0&I&0 \\
			0&-V_{1,2}^{(1)}{}^T&0&I
		\end{bmatrix}}_\text{$\widehat{K}$}
		\underbrace{
			\begin{bmatrix}
				x_1^{(1)}\\
				x_2^{(1)}\\
				y_1^{(1)} \\
				y_2^{(1)}
		\end{bmatrix}}_\text{$\widehat{x}$}
		=
		\underbrace{
			\begin{bmatrix}
				F_1\\
				F_2\\
				0\\
				0
		\end{bmatrix}}_\text{$\widehat{F}$}
\end{equation}
We now factorize the top diagonal block of $\widehat{K}$ which consists of $K_1^{(1)}$ and $K_2^{(1)}$. Since this subblock of $\widehat{K}$ is a block diagonal matrix, this means that we only need to factorize $K_1^{(1)}$ and $K_2^{(1)}$. After eliminating the top off diagonal block, we are left with the Schur complement:
\begin{equation}
	\label{eq:Schur}
	{S^{(1)}}=
	\begin{bmatrix}
		I & V_{2,1}^{(1)}{}^T (K_1^{(1)})^{-1} U_1^{(1)} \\
		V_{1,2}^{(1)}{}^T (K_2^{(1)})^{-1} U_2^{(1)} & I
	\end{bmatrix}
\end{equation}
All we have to do now, is to solve the Schur complement:
\begin{equation}
	\label{eq:y}
	S^{(1)}
	\begin{bmatrix}
		y_1^{(1)} \\
		y_2^{(1)}
	\end{bmatrix}=
	\begin{bmatrix}
		{V_{2,1}^{(1)}}^T  (K_1^{(1)})^{-1} F_1 \\ 
		{V_{1,2}^{(1)}}^T  (K_2^{(1)})^{-1} F_2
	\end{bmatrix}
\end{equation}
At this point, we can write $x_1^{(1)}$ and $x_2^{(1)}$ in terms of $(K_1^{(1)})^{-1}$ and $(K_2^{(1)})^{-1}$:
\begin{equation}
	\begin{bmatrix}
	x_1^{(1)}\\
	x_2^{(1)}
	\end{bmatrix}
	=
	\begin{bmatrix}
	(K_1^{(1)})^{-1} & 0\\
	0 & (K_2^{(1)})^{-1}
	\end{bmatrix}
	\begin{bmatrix}
	 	F_1 - U_1^{(1)} y_2^{(1)} \\
	 	F_2 - U_2^{(1)} y_1^{(1)}
	\end{bmatrix}
\end{equation}

Since, both ${K_1^{(1)}}$ and ${K_2^{(1)}}$ are HODLR matrices, we can apply the same procedure for factorizing them. Thus, we have arrived at a recursive algorithm for solving~\eqref{eq:modifiedSystem}. The factorization step corresponds to the computation and storage of all the terms that are independent of the right hand side (i.e., the Schur complements at all levels).  


\subsection{Algorithm Summary}
\label{sec:algorithmSummary}
	We now summarize the recursive HODLR direct solver algorithm. For a matrix such as $K\in \mathbb{R}^{n\times n}$, we have to carry out the following procedure at each recursion level $(p)$ for all $1\le i \le 2^p$:
\subsubsection{Factorize}
\begin{enumerate}
	\item Find the low-rank approximation of the off-diagonal blocks ($U_{2i-1}^{(p)}$, $U_{2i}^{(p)}$, $V_{2i-1,2i}^{(p)}$, $V_{2i,2i-1}^{(p)}$).
	\item Define $Z_1^0 = 0$. For each level $p$, starting at the top level ($p=0$), let:
		\begin{align}
		\begin{bmatrix}
		Z_{2i-1}^{(p+1)}\\
		Z_{2i}^{(p+1)}
		\end{bmatrix}=
		\left[
		\begin{array}{cc}
		U_{2i-1}^{(p+1)}&\multirow{2}{*}{$Z_{i}^{(p)}$}\\
		U_{2i}^{(p+1)}&
		\end{array}\right]
		\end{align}
	In the equation above, on the right-hand side, we are vertically concatenating two matrices to form a matrix at level $p+1$.
	\item Recursively solve the following equations:
		\begin{align}
			\label{eq:topDiag}
			\begin{bmatrix}d_{2i-1}^{(p+1)}&c_{2i-1}^{(p+1)} \end{bmatrix}&= (K_{2i-1}^{(p+1)})^{-1}Z_{2i-1}^{(p+1)}\\  \label{eq:bottDiag}
			\begin{bmatrix}d_{2i}^{(p+1)}&c_{2i}^{(p+1)} \end{bmatrix}&= (K_{2i}^{(p+1)})^{-1} Z_{2i}^{(p+1)}
		\end{align}
	where $d^{(p+1)}$ and $c^{(p+1)}$ correspond to the $U^{(p+1)}$ and $Z^{(p)}$ portion of the right hand sides respectively.
	\item Obtain $S_{i}^{(p)}$, using Equations~\eqref{eq:Schur} and~\eqref{eq:y}:
			\begin{equation}
				\label{eq:SchurSum}
				S_{i}^{(p)} =
				\begin{bmatrix}
					I & (V_{2i,2i-1}^{(p+1)})^T {d_{2i-1}^{(p+1)}}\\
					(V_{2i-1,2i}^{(p+1)})^T {d_{2i}^{(p+1)}}&I\\
				\end{bmatrix}
			\end{equation}
		\item Obtain $d_{i}^{(p)}$, $c_{i}^{(p)}$ for $p\ge1$ using:
		\begin{equation}
			\begin{bmatrix}
				d_{i}^{(p)} &c_{i}^{(p)}
			\end{bmatrix}=\left(
			  I-
			  \begin{bmatrix}
			  	0 & d_{2i-1}^{(p+1)}\\
				d_{2i}^{(p+1)}&0
			  \end{bmatrix}
			  (S_i^{(p)})^{-1}
			    \begin{bmatrix}
			  	{V_{2i,2i-1}^{(p+1)}}^T & 0\\
				0 & {V_{2i-1,2i}^{(p+1)}}^T
			  \end{bmatrix}\right)
			  \begin{bmatrix}
			  	c_{2i-1}^{(p+1)}\\
				c_{2i}^{(p+1)}
			  \end{bmatrix}
		\end{equation}
\end{enumerate}
	
\subsubsection{Solve}
\begin{enumerate}
	\item Define $z_1^0 = F$. For each level $p$, starting at the top level ($p=0$), let:
		\begin{align}
		\begin{bmatrix}
		z_{2i-1}^{(p+1)}\\
		z_{2i}^{(p+1)}
		\end{bmatrix}=z_i^{p}
		\end{align}
	\item Recursively solve the following equations:
		\begin{align}
			x_{2i-1}^{(p+1)}&= (K_{2i-1}^{(p+1)})^{-1} z_{2i-1}^{(p+1)}\\
			x_{2i}^{(p+1)}&= (K_{2i}^{(p+1)})^{-1} z_{2i}^{(p+1)}
		\end{align}
	\item Obtain $x_{i}^{(p)}$ for $p\ge0$ using:
		\begin{equation}
			x_{i}^{(p)}=\left(I-
			  \begin{bmatrix}
			  	0 & d_{2i-1}^{(p+1)}\\
				d_{2i}^{(p+1)}&0
			  \end{bmatrix}
			  (S_i^{(p)})^{-1}
			    \begin{bmatrix}
			  	{V_{2i,2i-1}^{(p+1)}}^T & 0\\
				0 & {V_{2i-1,2i}^{(p+1)}}^T
			  \end{bmatrix}\right)
			  \begin{bmatrix}
			  	x_{2i-1}^{(p+1)}\\
				x_{2i}^{(p+1)}
			  \end{bmatrix}
		\end{equation}
		Note that $(S_i^{(p)})^{-1}$ was previously computed and this step is therefore only a series of matrix-matrix products. Hence, the computational cost is small compared to the previous factorization.
\end{enumerate}

\subsection{Solver Computational Cost}
\label{sec:cost}
Assuming we use a fast ($\mathcal{O}(n)$) low-rank approximation scheme, the cost of constructing and storing an HODLR matrix is $\mathcal{O}(nr\log(n))$~\cite{SivaFDS}, where $r$ is the rank of approximation. 
Looking at the procedure described in Section~\ref{sec:algorithmSummary}, we can write the following:
\begin{equation}
	\label{eq:cost}
	C^{(p)}(r,s,n) = 2 C^{(p+1)}\left(r,s+r,\frac{n}{2} \right)
	+ \mathcal{O}(nr^2) + \mathcal{O}(nsr)
\end{equation}
where $C^{(p)}(r,s,n)$ is the computational cost associated with solving an $n\times n$ HODLR matrix at level $p$ with $s$ right hand sides and off-diagonal blocks of rank $r$. Equation~\eqref{eq:cost} suggests that the cost of solving a HODLR matrix at level $p$ with $s$ right hand sides is made up of three contributions. The first contribution is associated with solving the two diagonal blocks at the lower level $(p+1)$ with $s+r$ right hand sides. The second contribution comes from constructing the Schur complement ${S^{(p)}}$ (Equation~\eqref{eq:Schur}) and the third contribution is the cost of constructing the right hand side of Equation~\eqref{eq:y}. Writing Equation~\eqref{eq:cost} as a sum, we have:
\begin{equation}
	\label{eq:costResult}
	  C^{(0)}(r,s,n) = \sum_{p=1}^{\log \left( \frac{n}{r} \right )} \mathcal{O}(pnr^2+nsr) 
\end{equation}

If the off-diagonal rank is constant throughout various levels in the HODLR tree, the computational cost of the algorithm is $\mathcal{O}(r^2n\log^2(n))$ according to Equation~\eqref{eq:costResult}. 

However, in many practical cases, the rank decays from root to leaves in the HODLR tree. Assume we can approximate $r$ as $O(n_p^{{1/2}})$ where $n_p$ is the size of a block at level $p$. Then, we have: $r_p = O(\frac{r_1} { 2^{{p/2}}})$, where $r_1$ is the rank at the top level.
According to Equation~\eqref{eq:cost}, the total computational cost involves two sums:
\begin{gather*}
\sum_{p=1}^{\log \left( \frac{n}{r} \right )} r_p^2 = O(r_1^2) \\
\sum_{p=1}^{\log \left( \frac{n}{r} \right )} r_p 
\sum_{q=2}^p (s + r_q) = O(\sum_{p=1}^{\log \left( \frac{n}{r} \right )} r_p (s+r_1))
= O( r_1 (s+r_1) )
\end{gather*}
Note in particular that the second sum is $O(r_1^2)$ instead of $O(r^2 \log^2 n)$. Finally:
\begin{equation}
	\label{eq:costResultr}
	  C^{(0)}(r,s,n) = \mathcal{O}(nr^2) 
\end{equation}
This result shows that in cases where the off-diagonal rank is decreasing, HODLR solvers can become very efficient and can compete with HSS solvers.

 
\section{Low-Rank Approximation Schemes}
\label {sec:introLR}
In this section, we discuss the various low-rank approximations schemes used for obtaining a low-rank representation of the off-diagonal blocks of the HODLR matrices in consideration. Although a variety of low-rank approximation algorithms (SVD, rank revealing LU, rank revealing QR, randomized algorithms, etc) are available, we require a scheme that has a computational cost of $\mathcal{O}(rn)$ where $r$ is the rank of approximation and $n$ is the size of the matrix. In the context of this work, we cannot use randomized SVD methods since no fast matrix-vector product algorithm applies in our benchmark settings. This limits our choices to methods like Chebyshev, partial pivoting ACA (Section~\ref{sec:ACA}) and the pseudo-skeleton low-rank approximation algorithm (Section~\ref{sec:pseudo_boundary}). Each of these methods has certain drawbacks:
\begin{itemize}
 \item The Chebyshev low-rank approximation algorithm is only suited to cases dealing with interaction of points via smooth kernels. 
 \item The partial pivoting ACA algorithm works well when the leverage score of the matrix~\cite{leverageScore} is uniform. That is, all rows and columns have fairly the same importance when constructing the low-rank approximation.  However, in cases where certain rows or columns play a special role and are critical to include in the low-rank approximation, ACA might fail to properly identify them, resulting in an inaccurate low-rank approximation.
 \item The accuracy of the pseudo-skeleton low-rank approximation scheme strongly depends on the method used for selecting rows and columns.
\end{itemize}
In order to construct a fast and robust low-rank approximation scheme, we introduce a method for selecting rows and columns in the pseudo-skeleton low-rank approximation algorithm. We call this new method the boundary distance low-rank approximation scheme (BDLR). 


\subsection{ACA Low-Rank Approximation}
\label{sec:ACA}
We use the ACA algorithm with partial pivoting as described by Rjasanow~\cite{ACA}. This algorithm is an algebraic low-rank approximation scheme and works on any dense matrix without any prior knowledge of the matrix. Both full pivoting and partial pivoting ACA search the matrix or the remaining Schur complement for the largest entry and use this entry as the pivot. The full pivoting algorithm, similar to rank revealing LU, scans all the matrix entries. Partial pivoting ACA avoids this expensive search by looking at the largest entry in a single row/column at each step.
The partial pivoting ACA algorithm has a cost of $\mathcal{O}(r(m+n))$, for a matrix $A\in \mathbb{R}^{m\times n}$~\cite{ACA}, where $r$ is the rank of approximation.


\subsection{Randomized Algorithms}
\label{sec:Random}
Randomized algorithms as described by~\cite{RndSummary,Rnd1,Rnd2,Rnd3} arrive at a low-rank approximation of matrix $A$ by forming a lower dimensional matrix $Y$ obtained from sampling rows and/or columns of the original matrix or by applying random projections to matrix $A$. They then obtain the orthonormal basis $Q$ for the range of $Y$ and approximate $A$ as:
\begin{equation}
	A\approx QQ^TA
\end{equation}
For a matrix of size $n\times n$, and without a fast matrix-vector product, these methods have a computational cost of $\mathcal{O}(n^2)$. Otherwise, the cost can be brought down to $\mathcal{O}(n)$ or $\mathcal{O}(n\log n)$. 

\subsection{Pseudo-Skeleton and Boundary Distance Low-Rank Approximation}
\label{sec:pseudo_boundary}
In order to construct a fast and accurate solver, we need an accurate and robust method to construct low-rank approximations. 
As we will show, BDLR is very robust and leads to accurate low-rank approximations. It works well in problems where the matrix can be related to a Green's function. (This is true for all linear PDE problems. Note that the Green's function needs to be smooth, with a singularity at the origin). In that case, large entries correspond to points close in space, which we associate as a simplification to nodes in the graph that are connected by few edges. Although this is a simple heuristic, it worked very well in our examples and allowed us to efficiently form accurate low-rank approximations.

The BDLR algorithm is a row  and column selection algorithm in the pseudo-skeleton low-rank approximation scheme. The pseudo-skeleton algorithm allows us to construct a low-rank approximation of a matrix by choosing a subset of rows and columns of that matrix. As mentioned in~\cite{pseudoSkeleton}, for a low-rank matrix $A$, if we pick a set of row indices ($i \in I=\{i_1,...,i_r\}$) and a set of column indices ($j \in J=\{j_1,...,j_r\}$) and define matrices $C$ and $R$ such that :
\begin{align}
R &=  A(I,:)\\
C &=  A(:,J)
\end{align}
Then, we can approximate $A$ to be :
\begin{equation}
A\approx C{\widehat{A}}^{-1}R
\end{equation}
where ${\widehat{A}} = A(I,J)$. If ${\widehat{A}}$ is not a square matrix or rank deficient, the Moore-Penrose pseudoinverse is needed for ${\widehat{A}}^{-1}$. In order to achieve a certain accuracy, one can increase the number of chosen rows and columns until the desired accuracy is reached. To monitor the error in the scheme, we pick rows and columns that are not in the set of rows and columns already chosen for low-rank approximation. We then monitor the relative Frobenius norm error on these rows and columns and increase the rank of the approximation until the relative Frobenius norm error falls below a certain tolerance.

For a rank $r$ pseudo-skeleton low-rank approximation, the inversion of ${\widehat{A}}$ has a computational cost of $\mathcal{O}(r^3)$. Monitoring the error has a computational cost of $\mathcal{O}(mr+nr-r^2)$ for $A\in \mathbb{R}^{m\times n}$. Thus, this method has an asymptotic complexity of $\mathcal{O}(nr)$.

As mentioned in Section~\ref{sec:intro}, we are predominantly interested in solving dense frontal matrices arising from the multifrontal elimination process of sparse finite-element matrices. In this case, every frontal matrix has a corresponding sparse matrix, which is a diagonal sub-block of the original finite-element matrix. This sparse matrix describes a graph that has rows and columns of the dense matrix as its vertices and the edges in this graph correspond to nonzero entries in the sparse matrix and describe the connection between these points.  We use this graph in constructing the low-rank approximation of the off-diagonal blocks.

Entries in dense matrix blocks that correspond to FEM or BEM applications can be related to the inverse of a Green's function. The Green's function is large at short distances and then decays smoothly. We have a similar behavior for our dense blocks. Hence, we want to identify row/column pairs corresponding to large entries. These correspond to nodes in the graph that are close, that is connected by few edges. Therefore we use the distance between a row vertex in the graph and the column vertex set (e.g., if the vertex corresponds to a row, we consider the distance to the set of vertices associated with the columns, and vice versa) as a good criterion to determine whether to pick a row/column or not. 

\begin{figure}[htbp]
	\centering
	\subfigure[Full Matrix Representation]{
	\includegraphics[scale=1]{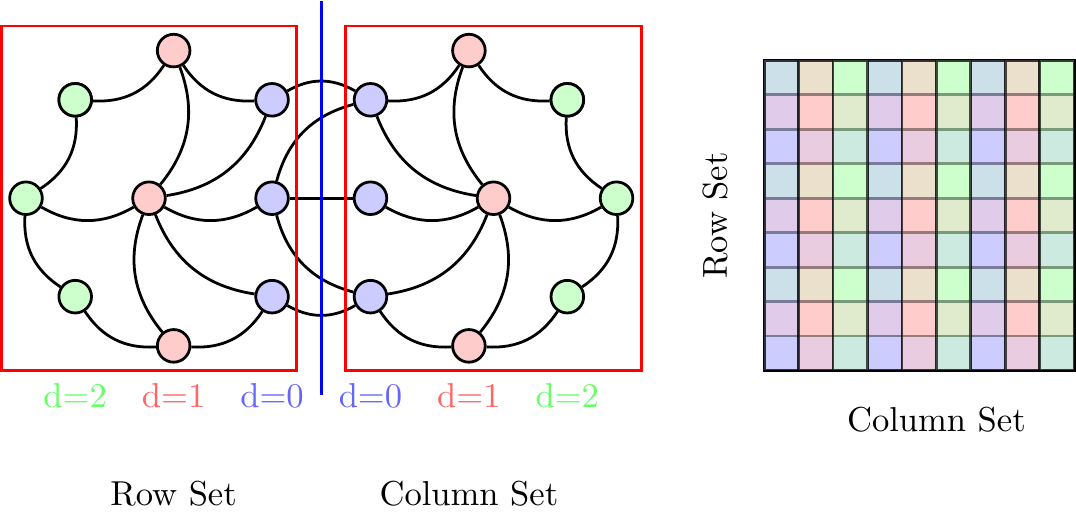}
	\label{fig:BDLR_Full}
	}
	\subfigure[Low-Rank Matrix Representation]{
	\includegraphics[scale=1]{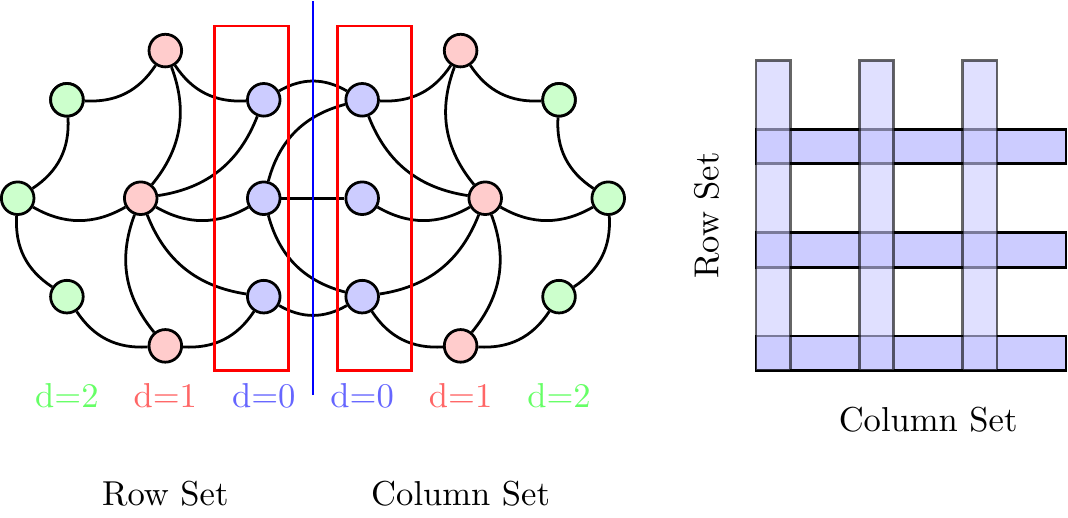}
	\label{fig:BDLR_LR}
	}
	\caption{Classification of vertices based on distance from the other set.}
\end{figure}

For a set of row (column) vertices, we define the boundary vertices as the subset of vertices for which there exists an edge in the interaction graph connecting them to a vertex in the column (row) set. Figure~\ref{fig:BDLR_Full} shows an example of a matrix which corresponds to the interactions of a set of row points with a set of column points. In this particular example, the blue vertices are the boundary vertices. That is, they are the vertices closest to the boundary between the row and column set of points.

Now that we have defined the boundary nodes, we can designate an index $d$ for every vertex in the row (column) set. This index is defined as the distance of a vertex to the vertices in the boundary set. In order to construct the low-rank approximation, we choose rows and columns based on their $d$ index value. That is, we first choose rows (columns) that are in the boundary set ($d = 0$). We then add rows (columns) with a distance of one to the boundary ($d = 1$). For example, in Figure~\ref{fig:BDLR_Full}, the green nodes are labeled ($d = 1$) as they are separated from the blue boundary nodes ($d = 0$) with only one edge. We continue adding points based on the $d$ index, until we reach the desired accuracy.  Figure~\ref{fig:BDLR_LR} shows that the BDLR algorithm approximates the interaction of a set of row and column nodes with the interaction of the ones that are closest to the boundary (interaction of blue nodes).

As mentioned above, calculating the pseudo skeleton low-rank approximation requires us to calculate the pseudoinverse of ${\widehat{A}}$. For the BDLR algorithm, instead of using the SVD for calculating the pseudoinverse (${\widehat{A}}^{-1}$), we use a full pivoting LU factorization, which is slightly cheaper:
\begin{equation}
	\label{eq:Ahat}
	{\widehat{A}} = P^{-1}LUQ^{-1}
\end{equation}
where $P$ and $Q$ are permutation matrices. Let r be the rank of $\widehat{A}$. Define $\widetilde{R}$ and $\widetilde{C}$ as:
\begin{align}
  \widetilde{C}& =  (CQ)(:,1:r) (U(1:r,1:r))^{-1}\\
  \widetilde{R}& =  (L(1:r,1:r))^{-1} (PR)(1:r,:)
\end{align}
where $C$ and $R$ are the subset of columns and rows we have picked using the BDLR scheme. We then have:
\begin{equation}
	A\approx \widetilde{C} \; \widetilde{R}
\end{equation}
$(U(1:r,1:r))^{-1}$ and $(L(1:r,1:r))^{-1}$ correspond to lower-triangular solves. The inverse matrices are not explicitly computed.

\section{Application for Multifrontal Solve Process}
\label{sec:MF}

In this section, we demonstrate how our fast dense solver algorithm can be applied to a sparse multifrontal solve process. We will not explain the multifrontal algorithm in detail. For a detailed review of the multifrontal method see~\cite{MFReview}. We applied our fast solver as described in Section~\ref{sec:directIterative} to a variety of 3D finite-element problems. 
We investigate frontal matrices at various levels of the sparse matrix elimination tree corresponding to the elasticity equation. 
We use SCOTCH~\cite{SCOTCH} to do the reordering in the sparse multifrontal solver. Our goal is to apply our fast dense solver to the dense frontal matrices obtained in the multifrontal elimination process of a sparse finite-element matrix, and speed up the multifrontal algorithm to approximately $\mathcal{O}(N^{4/3})$. The results shown in this paper can be viewed as a proof of concept of this idea. We should also mention that the approach presented in this article is fully general. We use SCOTCH~\cite{SCOTCH}, (which can partition any graph) to obtain the separators and the resulting separators can always be handled by our algorithm, without any change.

\subsection{Elasticity Equation for a 3D Beam and a Cylinder Head Geometry}
\begin{figure}[htbp]
	\centering
	\subfigure[Beam]{
	\includegraphics[width=200pt,height=180pt]{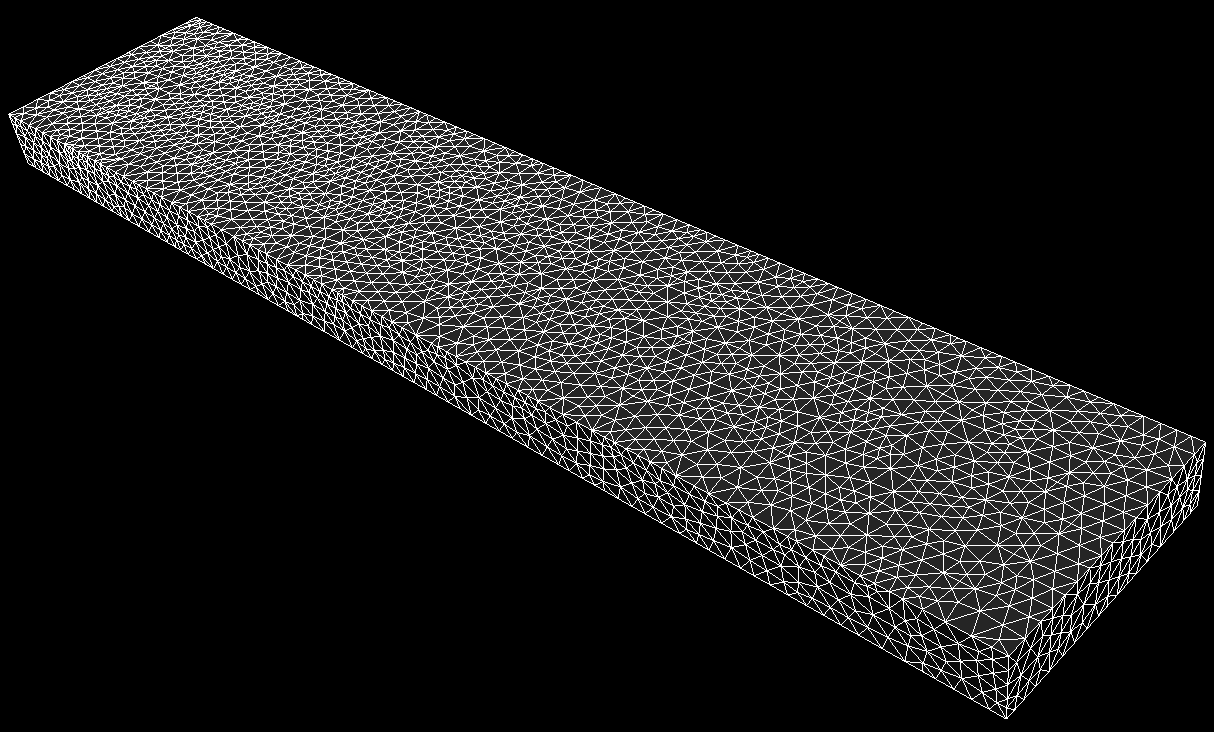}
	\label{fig:beam}}
	\subfigure[Cylinder Head]{
	\includegraphics[width=200pt,height=180pt]{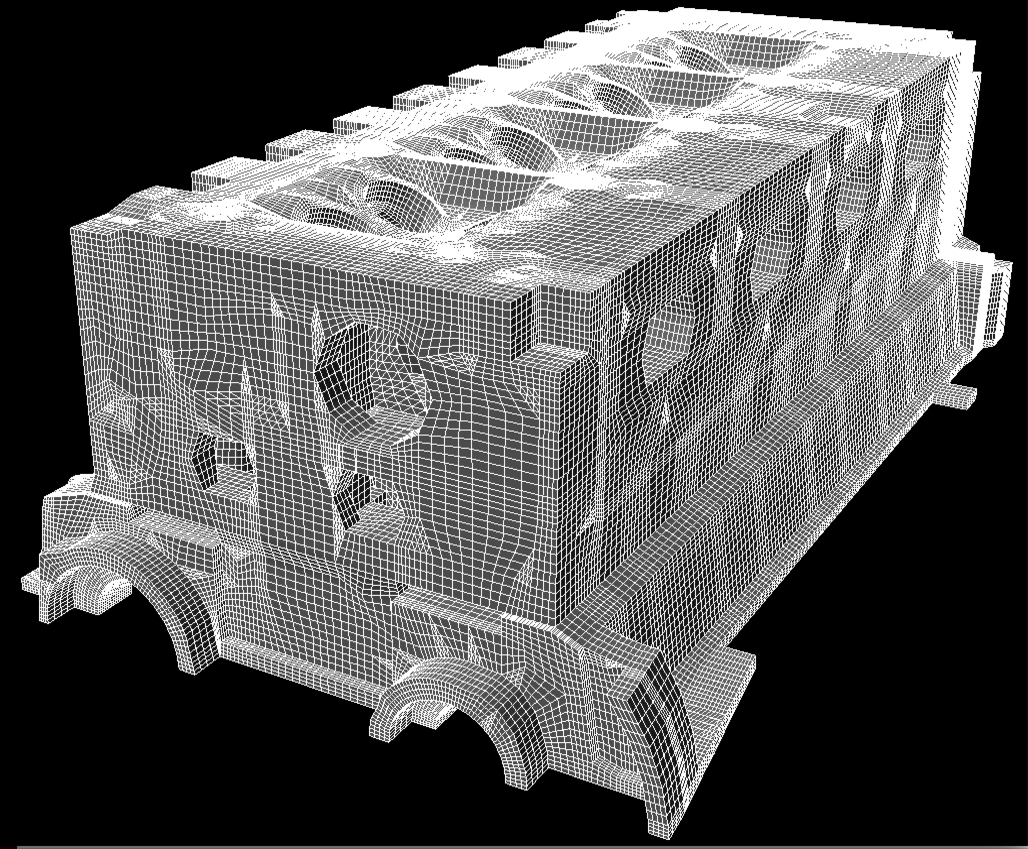}
	\label{fig:cylinderHead}}
	\caption{3D unstructured mesh for the beam and cylinder head geometries. .}
\end{figure}

We consider the 3D Navier-Cauchy elastostatics equations with a beam geometry (figure~\ref{fig:beam}):
\begin{equation}
	(\lambda + \mu)\nabla(\nabla \cdot \boldsymbol{u})+\mu \nabla^2\boldsymbol{u}+\boldsymbol{F} = 0
	\label{eq:NavierCauchy}
\end{equation}
where $\boldsymbol{u}$ is the displacement vector and $\lambda$ and $\mu$ are Lam\'{e} parameters. For the beam geometry, we use 10-node tetrahedral elements (see for example Section~ 10.2 of this document\footnote{\url{http://www.colorado.edu/engineering/CAS/courses.d/AFEM.d/AFEM.Ch10.d/AFEM.Ch10.pdf}}) to discretize the above equation. For the cylinder head geometry, the mesh is composed of 8-node hexahedral, 6-node pentahedral and 4-node tetrahedral solid elements, and also 3-node shell elements. Figures~\ref{fig:beam} and~\ref{fig:cylinderHead} show a sample beam and cylinder head geometry respectively. As can be seen, the meshes are unstructured for both geometries.

\subsection{FETI-DP Solver for a 3D Elasticity Problem}

Domain decomposition (DD) methods solve a problem by splitting it into several subdomains. Local problems are solved on each subdomain and a global linear system is used to couple these local solutions into a global solution for the entire problem~\cite{DD}. FETI methods~\cite{FETI_DP1,FETI_DP2} are a family of domain decomposition algorithms with Lagrange multipliers that have been developed for the fast sequential and parallel iterative solution of large-scale systems of equations arising from the finite-element discretization of partial differential equations~\cite{FETI_DP1}.

\begin{figure}[htbp]
	\centering
	\subfigure[Structured Mesh]{
	\includegraphics[width=200pt,height=180pt]{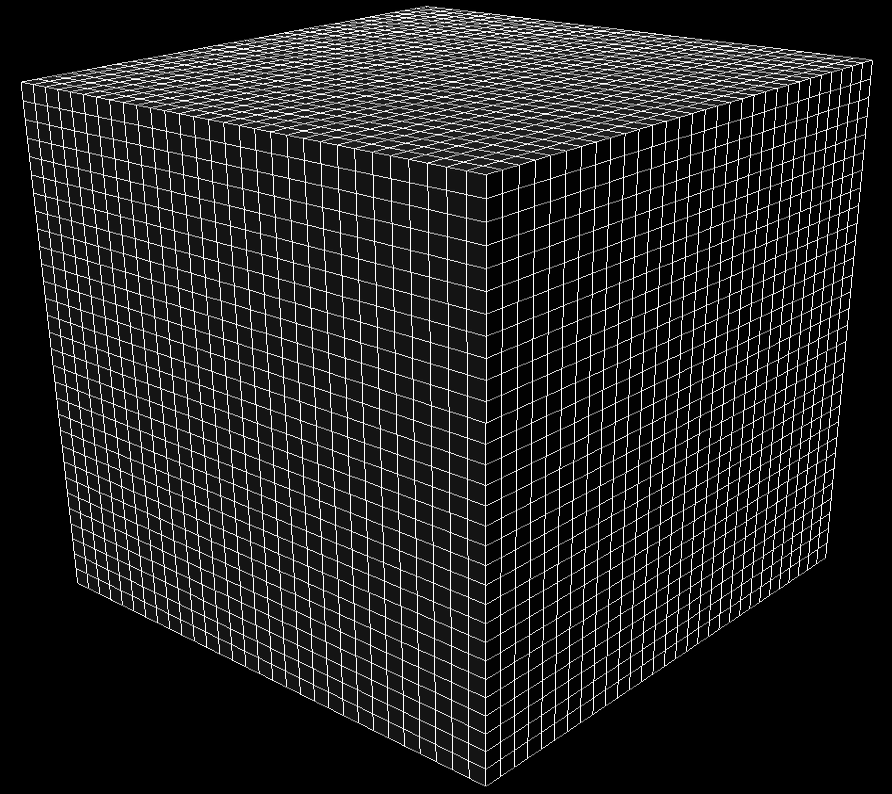}
	\label{fig:cube}
	}
	\subfigure[Unstructured Mesh]{
	\includegraphics[width=200pt,height=180pt]{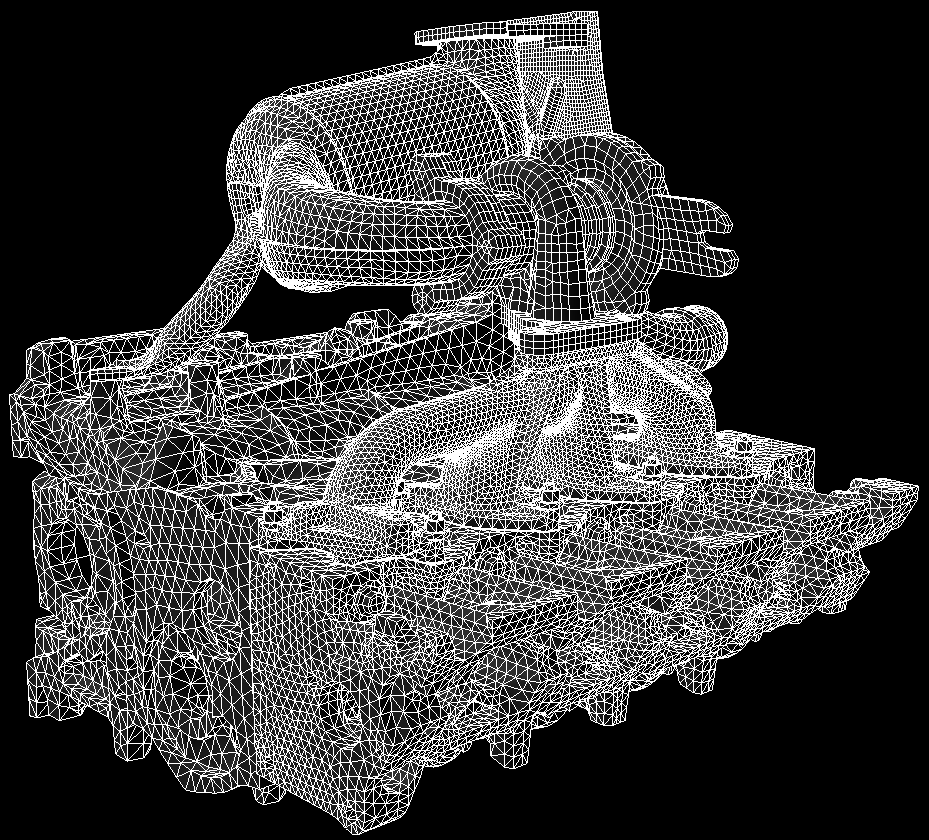}
	\label{fig:engine}
	}
	\caption{FETI-DP benchmark meshes. Figure (a) shows a structured and figure (b) shows an unstructured 3D FETI-DP mesh.}
\end{figure}

In this article, we consider two sparse local FETI-DP matrices arising from the finite-element discretization of an elasticity problem in three dimensions. The first matrix corresponds to solving the elasticity equation with a structured mesh in three dimensions (figure~\ref{fig:cube}) while the second matrix corresponds to solving the same problem using the geometry of an engine in an unstructured mesh (figure~\ref{fig:engine}). Both matrices correspond to the stiffness matrix of one subdomain of a linear elastic 3D solid finite element model (Equation~\eqref{eq:NavierCauchy}) of their respective geometry. The discretization for the cube geometry uses 8-node (trilinear) hexahedral elements (see for example Section 11.3 of this online document\footnote{\url{http://www.colorado.edu/engineering/CAS/courses.d/AFEM.d/AFEM.Ch11.d/AFEM.Ch11.pdf}}) while the discretization for the engine geometry uses 10-node tetrahedral elements (see for example Section 10.2 of this document\footnote{\url{http://www.colorado.edu/engineering/CAS/courses.d/AFEM.d/AFEM.Ch10.d/AFEM.Ch10.pdf}}).

\section{Numerical Benchmarks}
\label{sec:benchmarks}

In this section we show some numerical results and benchmarks of our code. As our code uses the Eigen C\texttt{++} library for matrix manipulations, we use the Eigen direct solvers as benchmark references.

\subsection{Elasticity Equation for a 3D Beam and a Cylinder Head Geometry}
We apply our solvers to frontal matrices arising from the multifrontal elimination of 3D elastostatics sparse matrices (Figures~\ref{fig:beam},~\ref{fig:cylinderHead}).
We compare the fast BDLR direct solver and the ACA direct solver as preconditioners to the GMRES iterative scheme. Because of the particular geometry of the beam mesh, all frontal matrices are relatively small ($\leq2$K) for this particular case.

\begin{figure}[htbp]
\centering
		\includegraphics{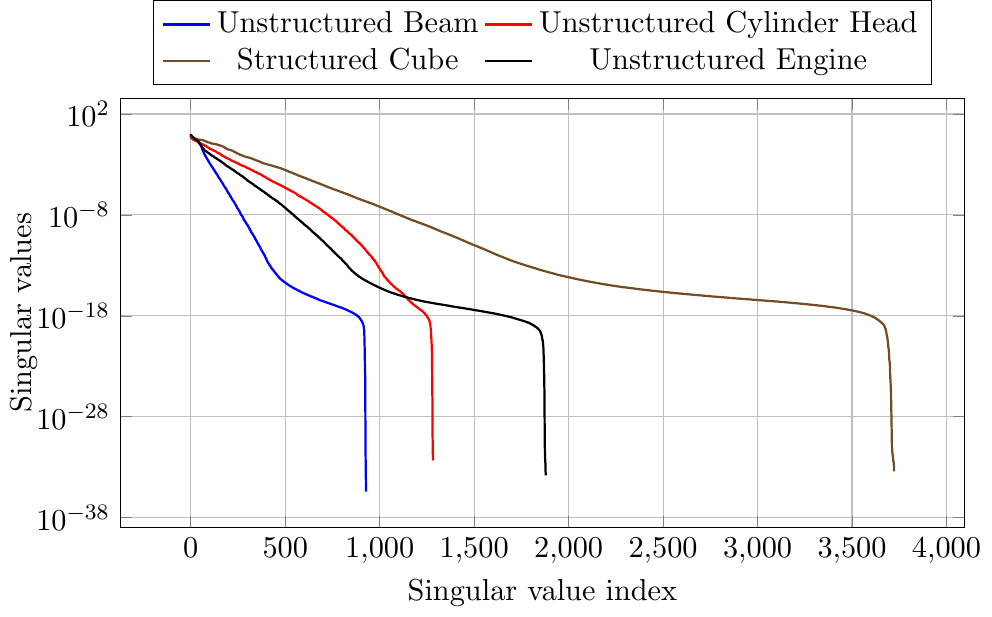}

\caption{Singular value decay for a variety of sample off-diagonal blocks of frontal matrices. The beam, cylinder head, cube and engine geometries correspond to blocks of size 0.95K, 1.3K, 3.75K and 1.9K respectively. }

\label{fig:SVD}
\label{fig:beamRank}
\end{figure}

As can be seen in Figure~\ref{fig:SVD}, the singular values of a sample frontal matrix off-diagonal block decay rapidly and the block is in fact low-rank. Figures~\ref{fig:beamRowDist} and~\ref{fig:beamColDist} show the distance of row (column) index of each pivot obtained in the full pivoting LU factorization from the boundary between the row and column sets of vertices in the interaction graph for the beam problem. As we expected, larger pivots correspond to rows and columns that are closer to the boundary. Figures~\ref{fig:beamPoints} and~\ref{fig:cylinderPoints} compares the relative error in approximating the top off-diagonal block using SVD versus the BDLR approximation for the beam and cylinder head geometry respectively. That is, each point (x,y) in this plot represents the relative error in approximation (y) if we wanted a rank (x) approximation using one of the low-rank approximation algorithms. Needless to say, this corresponds to choosing the top singular values in the SVD decomposition and choosing rows and columns that are closest to the boundary in the BDLR approximation. As can be seen in the plot, the curves associated with the BDLR scheme have a tolerance ($\epsilon$). This means that after the LU factorization of $\widehat{A}$ (see Section~\ref{sec:pseudo_boundary}), we only keep rows and columns corresponding to pivots that are larger than $\epsilon$ times the magnitude of the largest pivot. We use this convention for all BDLR approximations in this paper. We can observe that as we decrease $\epsilon$, we obtain a more accurate low-rank representation via the BDLR algorithm for the beam geometry. For the more complicated cylinder head geometry, we see that in order to obtain a good approximation for low values of $\epsilon$, more rows and columns need to be included in the low-rank approximation which corresponds to a higher depth parameter ($d$) in the BDLR scheme.

\begin{figure}[htbp]

\centering
\subfigure[Unstructured Beam ]{
		\includegraphics{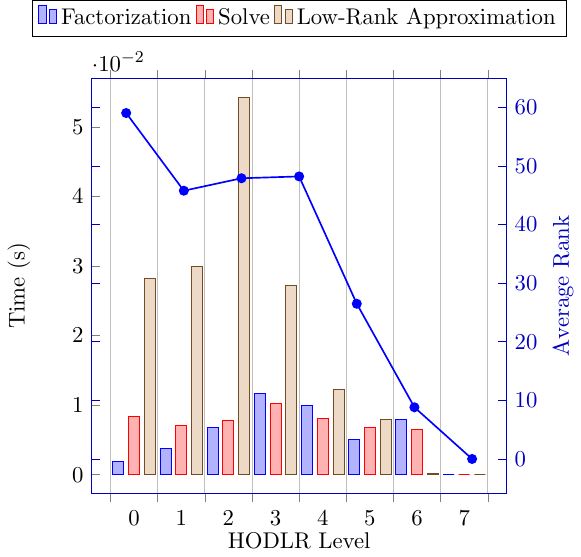}

	\label{fig:beamRankTiming}
}
\centering
\subfigure[Unstructured Cylinder Head]{
		\includegraphics{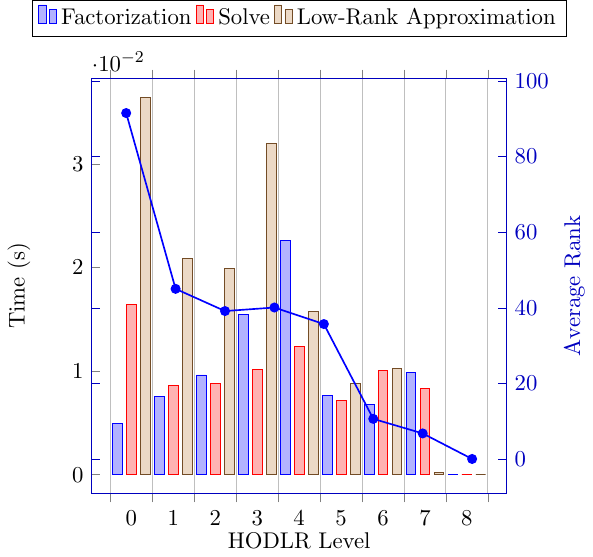}

	\label{fig:cylinderHeadRankTiming}
}
\centering
\subfigure[Structured Cube]{
		\includegraphics{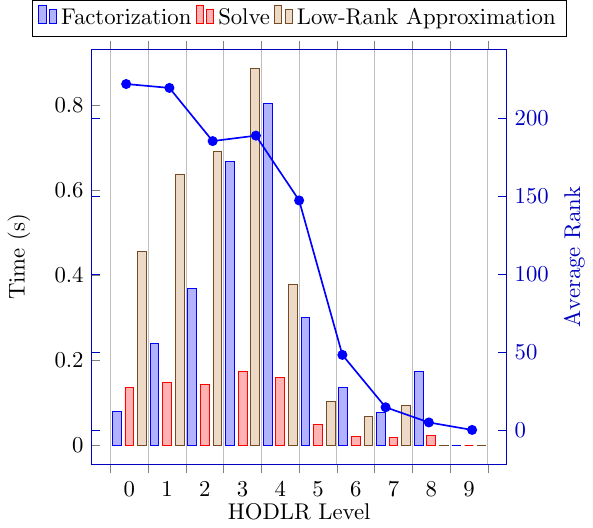}

	\label{fig:cubeRankTiming}
}
\centering
\subfigure[Unstructured Engine]{
		\includegraphics{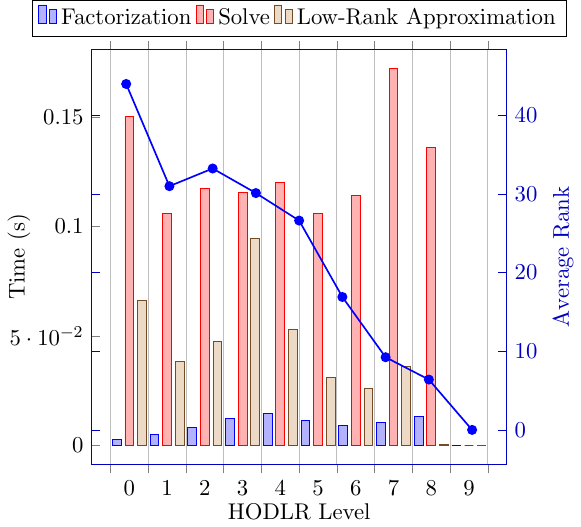}

	\label{engineRankTiming}
}

\caption{Off-diagonal rank and level by level timings for various frontal matrices. The off-diagonal ranks correspond to a BDLR low-rank approximation scheme with depth of $1$ and tolerance of $10^{-1}$. The left axis corresponds to runtimes for various stages in the solve process as a function of HODLR level. The right axis shows the off-diagonal rank versus HODLR level. a) A  frontal matrix corresponding to an unstructured beam geometry of size 1.9K. b) A frontal matrix corresponding to an unstructured cylinder head mesh of size 2.6K. c) A frontal matrix corresponding to a structured cube mesh of size 7.5K. d) A frontal matrix corresponding to an unstructured engine mesh of size 3.8K.}
\end{figure}

Figures~\ref{fig:beamRankTiming} and~\ref{fig:cylinderHeadRankTiming} show a level by level timing of the factorization, solve and low-rank approximation of the BDLR solver applied to sample frontal matrices corresponding to the beam and cylinder head geometries respectively. As can be seen, the off-diagonal rank decays from root to leaf which confirms our assumptions in Section~\ref{sec:cost}. Figures~\ref{fig:beamIter} and~\ref{fig:cylinderHeadIter} show a detailed convergence analysis and comparison between the BDLR and ACA solvers as preconditioners to the GMRES iterative scheme.

\subsection{FETI-DP Solver for a 3D Elasticity Problem}
\label{sec:FETIResults}
We apply the BDLR and ACA direct solver preconditioner to frontal matrices arising from the multifrontal elimination of local matrices in a FETI-DP solver. We considered two different classes of problems. One corresponds to solving the elasticity equation (Equation~\eqref{eq:NavierCauchy}) in a cube geometry with a structured mesh. The other corresponds to solving the same equation in an engine geometry with an unstructured mesh.

\begin{figure}[htbp]
\centering
\subfigure[Row Distance (Unstructured Beam)]{
		\includegraphics{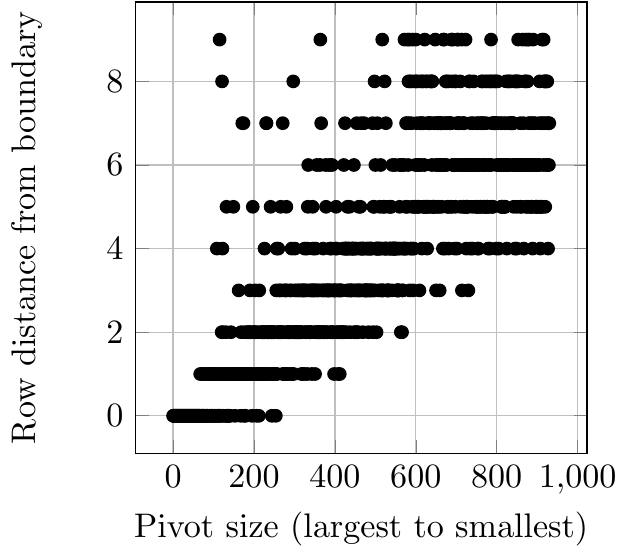}
	\label{fig:beamRowDist}
}
\subfigure[Col Distance (Unstructured Beam)]{
		\includegraphics{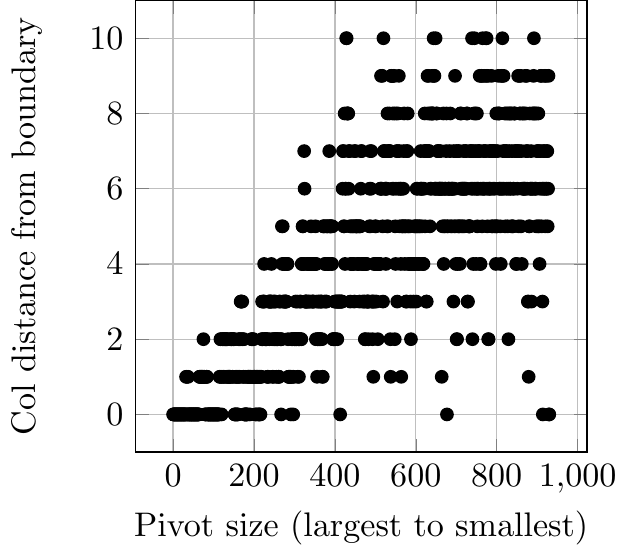}
	\label{fig:beamColDist}
}

\subfigure[Row Distance (Structured Cube)]{
		\includegraphics{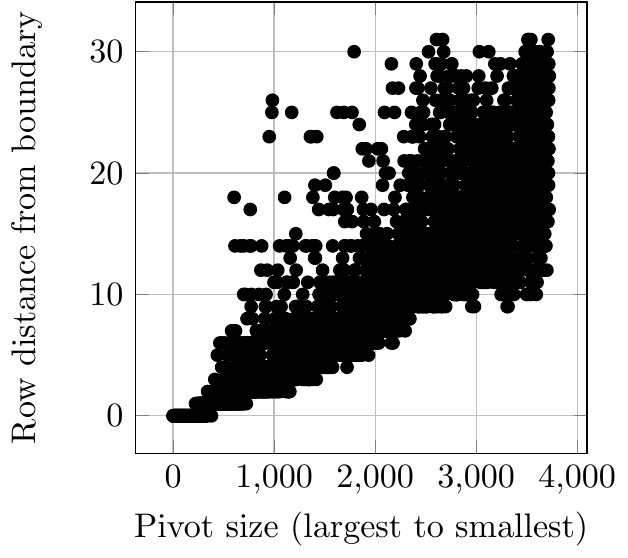}

	\label{fig:cubeRowDist}
}
\subfigure[Col Distance (Structured Cube)]{
		\includegraphics{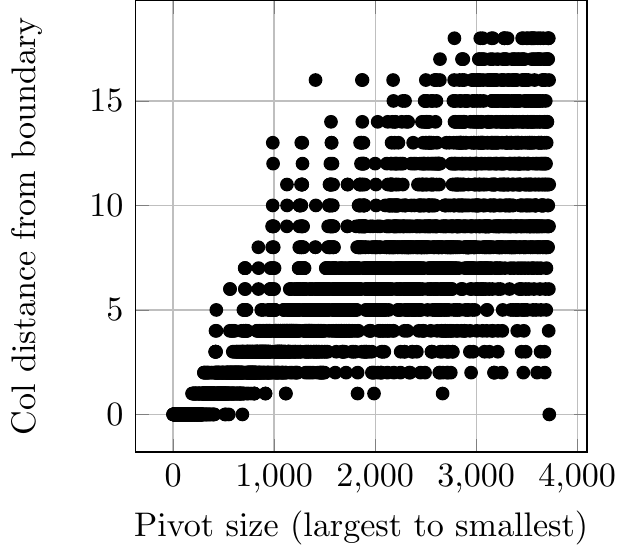}

	\label{fig:cubeColDist}
}

\subfigure[Row Distance (Unstructured Engine)]{
		\includegraphics{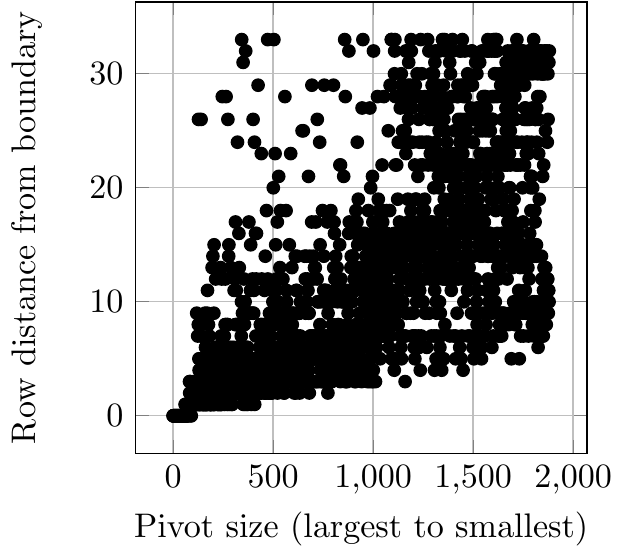}

	\label{fig:engineRowDist}
}
\subfigure[Col Distance (Unstructured Engine)]{
		\includegraphics{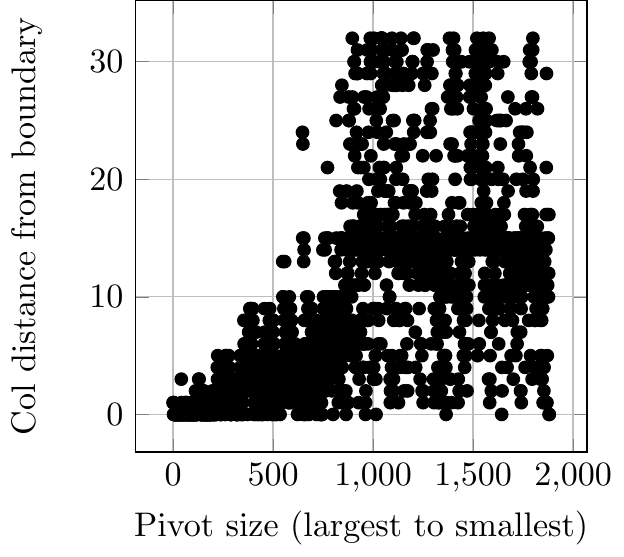}

	\label{fig:engineColDist}
}

\caption{Row (column) distance versus pivot size for a variety of off-diagonal blocks of sample frontal matrices. Row (column) distance is the distance corresponding to the row (column) index of a pivot from the boundary as defined in Figure~\ref{fig:BDLR_Full}. This graph shows that large pivots are near the boundary interface, whereas the pivot size decays as we move away. This justifies heuristically our approach with BDLR. a,b) An off diagonal block of an unstructured beam geometry frontal matrix of size 0.95K. c,d) An off diagonal block of an structured cube geometry frontal matrix of size 3.75K. e,f) An off diagonal block of an unstructured engine geometry frontal matrix of size 1.9K.}
\label{fig:cubeRank}
\end{figure}

Figures~\ref{fig:cubeRowDist} and ~\ref{fig:cubeColDist} show that the largest pivot values of a sample off-diagonal block of a frontal matrix arising from the cube geometry correspond to rows and columns that are closer to the boundary. Figures~\ref{fig:engineRowDist} and~\ref{fig:engineColDist} show that for the unstructured engine mesh, although most large pivots correspond to rows and columns near the boundary, there are some important rows and columns that are not included in the points closest to the boundary.

\begin{figure}[htbp]
\centering
\subfigure[Unstructured Beam]{
		\includegraphics{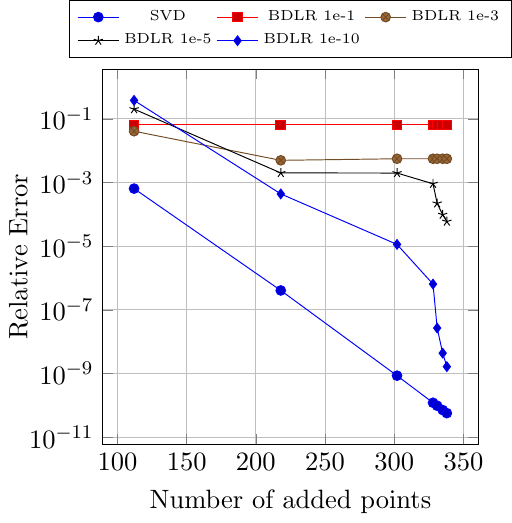}

	\label{fig:beamPoints}
}
\centering
\subfigure[Unstructured Cylinder Head]{
		\includegraphics{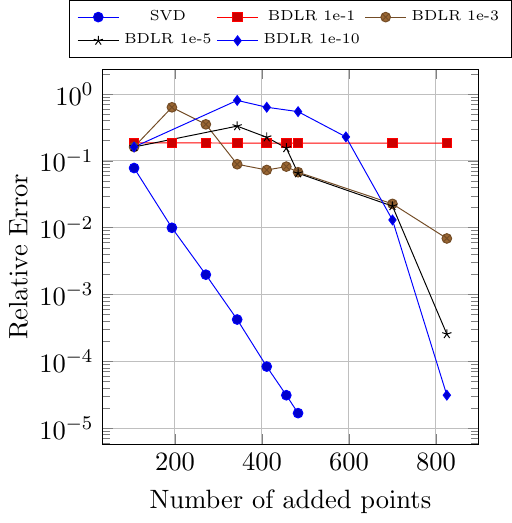}

	\label{fig:cylinderPoints}
}
\centering
\subfigure[Structured Cube]{
			\includegraphics{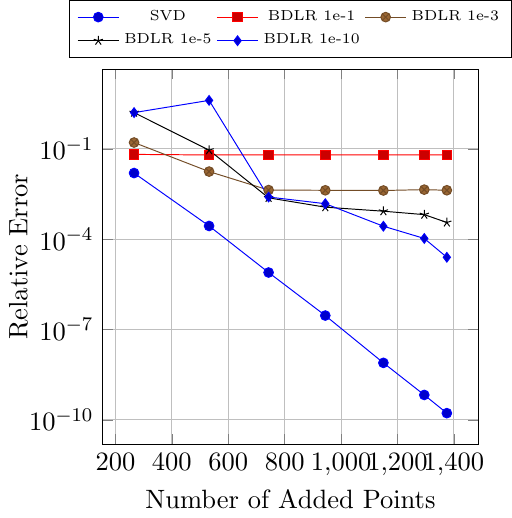}

	\label{fig:cubePoints}
}
\centering
\subfigure[Unstructured Engine]{
			\includegraphics{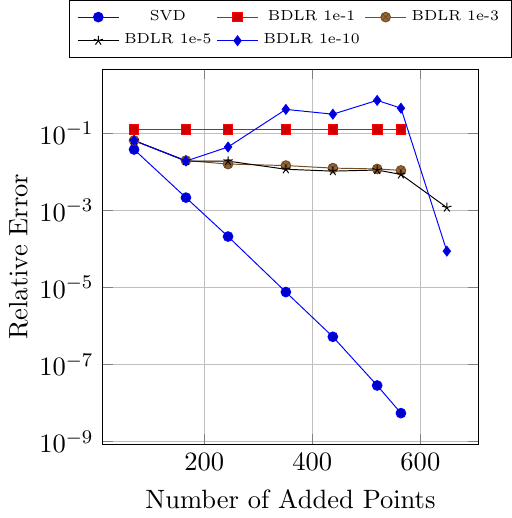}

\label{fig:enginePoints}

}

\caption{BDLR error vs SVD error for a variety of sample frontal matrix off-diagonal blocks. BDLR accuracy is used to truncate the pivots and number of points is the size of the block for full pivoting in the pseudo-skeleton approach (size of $\hat{A}$ in Equation~\eqref{eq:Ahat}). a) An off diagonal block of an unstructured beam geometry frontal matrix of size 0.95K. b) An off diagonal block of an unstructured cylinder head geometry frontal matrix of size 1.3K. c) An off diagonal block of an structured cube geometry frontal matrix of size 3.75K. d) An off diagonal block of an unstructured engine geometry frontal matrix of size 1.9K.}

\label{fig:engineRank}
\end{figure}

Figures~\ref{fig:cubePoints} show that the error in the BDLR method is comparable to the SVD (optimal) algorithm for the structured cube problem. Figure~\ref{fig:enginePoints} shows that similar to Figure~\ref{fig:cylinderPoints}, we need to include more points (rows and columns), in order to achieve an accurate low-rank approximation for $\epsilon = 10^{-10}$.
In other words, if there are insufficient rows and columns in the BDLR approximation, the matrix $\widehat{A}$ (see Section~\ref{sec:pseudo_boundary}) becomes low-rank and results in a LU factorization with vey small pivots. These small pivots are the cause of the large relative error as they become very large when inverted.

\begin{figure}[htbp]
\centering
\subfigure[Unstructured Beam]{
		\includegraphics{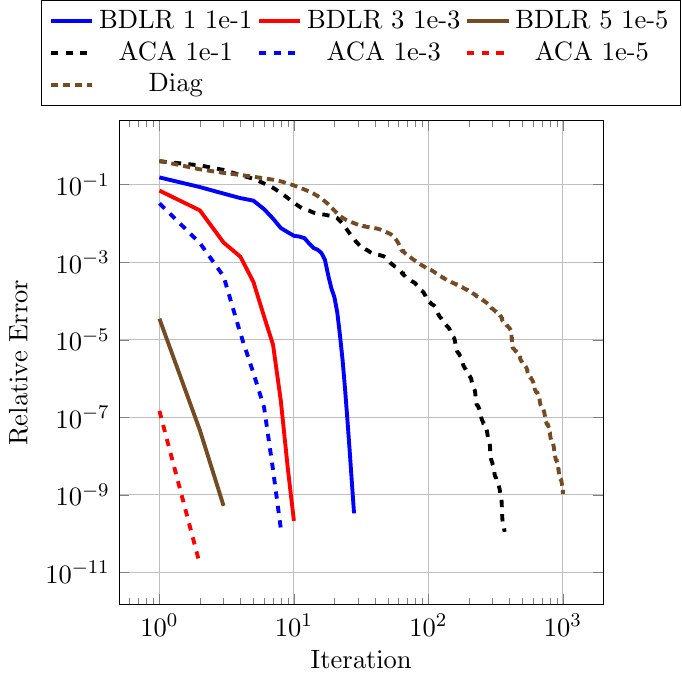}

	\label{fig:beamIter}
}
\centering
\subfigure[Unstructured Cylinder Head]{
		\includegraphics{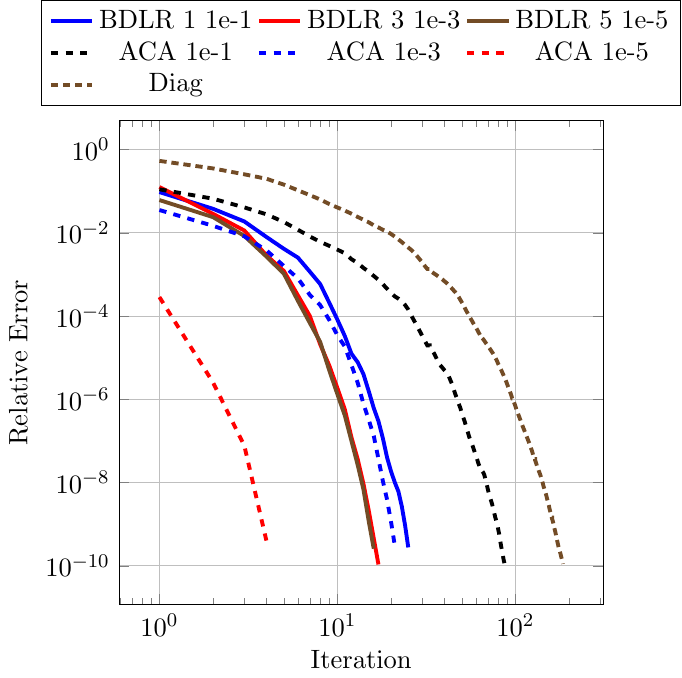}

	\label{fig:cylinderHeadIter}
}
\centering
\subfigure[Structured Cube]{
		\includegraphics{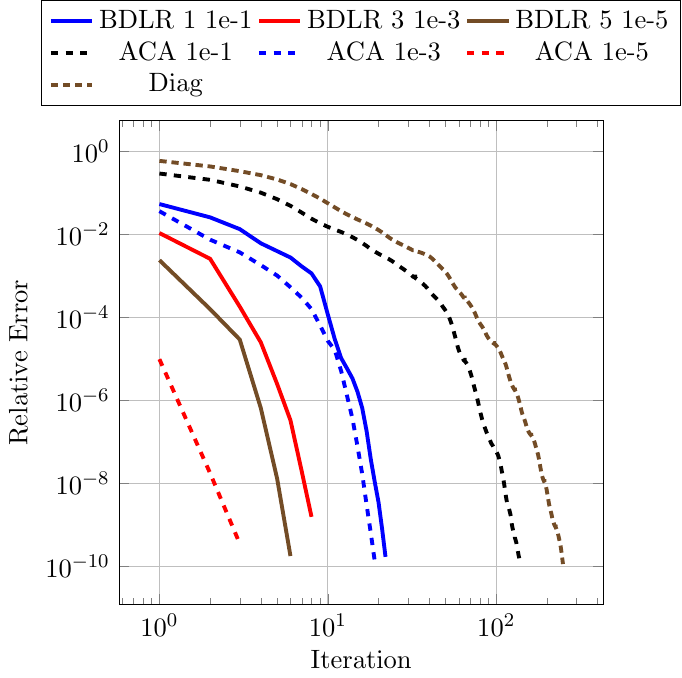}

	\label{fig:cubeIter}
}
\centering
\subfigure[Unstructured Engine]{
		\includegraphics{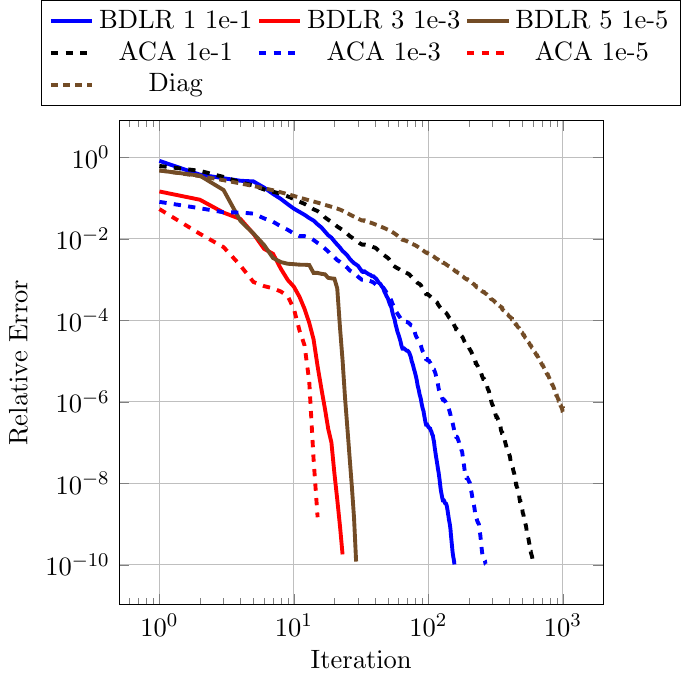}

	\label{fig:engineIter}
}

%

	\caption{Convergence analysis for BDLR and ACA preconditioners with the GMRES iterative scheme for a variety of frontal matrices. The curve labeled `Diag' corresponds to GMRES with diagonal preconditioning. a) A frontal matrix corresponding to an unstructured beam mesh of size 1.9K. b) A frontal matrix corresponding to an unstructured cylinder head mesh of size 2.6K. c) A frontal matrix corresponding to a structured cube mesh of size 7.5K. d) A frontal matrix corresponding to an unstructured engine mesh of size 3.8K.}
\label{fig:engineSpeed}
\end{figure}

Figures~\ref{fig:cubeIter} and~\ref{fig:engineIter} show the convergence rate of various BDLR and ACA direct solver preconditioners for a sample frontal matrix arising from the cube and engine mesh respectively. 

\subsection{Summary}
Table~\ref{table:summaryHODLRNum} summarizes the solver timings for various frontal matrices that we benchmarked. As can be seen, the iterative solve scheme with both a fast BDLR and ACA direct solver preconditioner can reach near machine accuracy much faster than a conventional LU solver in almost all cases. Furthermore, both BDLR and ACA achieve a relatively good speedup for all cases. However, for very large cases ($1.5$M structured cube and $2.3$M unstructured cylinder head), one can observe that  BDLR achieves higher speedup compared to ACA. One important point to note, is that convergence of both BDLR and ACA depends on the chosen parameters. For ACA, one can get better results by decreasing the tolerance. For BDLR, in order to achieve a given tolerance, one has to increase the depth parameter ($d$). It is possible for BDLR not to converge for a certain tolerance and a depth parameter. This is because the depth and accuracy are related. In particular, the efficiency of the method is sometimes found to degrade if we reduce $\epsilon$ too much without increasing $d$ sufficiently. This corresponds to the fact that we are trying to get a more accurate low-rank approximation but the pool of sample points is not sufficiently large to provide the desired accuracy. In that case, reducing $\epsilon$ may, in fact, lead to a degradation in the preconditioner, rather than an improvement. 



An important advantage of the BDLR algorithm is that the rows and columns required for constructing the low-rank approximation are known a priori based on the structure of the separator graph. As we will demonstrate in a future article, this will allow us to significantly accelerate the extend-add process and allows us to avoid constructing large dense frontal and update matrices as we will only keep track of rows and columns required by the BDLR algorithm.

%

\begin{table}[htbp]
	\scalebox{0.54}{
	\begin{tabular}{|c|c|c|c||c|c|c|c|c|c|c|c|c|c|c|c|c|c|c|c|}
	\hline
	 \multirow{2}{*}{Matrix} & \multirow{2}{*}{Mesh} & \multirow{3}{*}{Level} & \multicolumn{2}{c}{Matrix Size}\vline & \multicolumn{6}{c}{ACA}\vline & \multicolumn{6}{c}{BDLR} \vline&  \multirow{3}{*}{LU} & \multicolumn{2}{c}{Speed-up}\vline \\ \cline{4-17} \cline{19-20}
	 \multirow{2}{*}{Type} &  \multirow{2}{*}{Type}  & & \multirow{2}{*}{Sparse} & \multirow{2}{*}{Dense}& \multicolumn{2}{c}{1e-1}\vline & \multicolumn{2}{c}{1e-3}\vline & \multicolumn{2}{c}{1e-5}\vline & \multicolumn{2}{c}{1e-1}\vline & \multicolumn{2}{c}{1e-3}\vline & \multicolumn{2}{c}{1e-5}\vline &&\multirow{2}{*}{ACA}&\multirow{2}{*}{BDLR} \\ \cline{6-17}
	  &  & &  & & T & I & T & I & T & I &T&I&T&I&T&I&&& \\ \hline
	
	\multirow{11}{*}{FETI Local}&\multirow{10}{*}{Cube}
	   &1st&1.5M&23K& 1.32e2 & 223 & 2.85e2 & 72 & 7.82e2 &15&1.12e2&34&2.90e2&13&6.71e2&7&7.29e2&\bf5.52&\bf6.51\\ \cline{3-20}
          & &1st &\multirow{9}{*}{400K}&7.5K&6.99e0&141&1.78e1&20&4.29e1&3&8.32e0&23&2.28e1&9&4.77e1&7&2.38e1&\bf3.40&\bf2.86 \\ 
          & &2nd&&5.2K&2.29e0&77  &6.91e0&19&1.51e2&3&3.23e0&17&9.03e0&7&1.83e1&6&8.53e0&\bf3.72&\bf2.64 \\  
	  & &2nd&&5.0K&2.50e0&74&7.82e0&19&1.68e1&2&3.38e0&17&9.75e0&6&2.06e1&4&7.40e0&\bf2.96&\bf2.18\\
	  & &3rd&&2.0K&2.77e-1&45&6.15e-1&9&1.21e0&3&3.33e-1&12&7.86e-1&5&1.44e0&4&5.41e-1&\bf1.95&\bf1.62\\
	  & &3rd&&2.8K&4.34e-1&61&1.40e1 & 13 & 2.89e1 & 3 &7.17e-1&15&1.76e0&7&3.72e0&11&1.31e0&\bf3.01&\bf1.83 \\ 
          & &3rd&&2.2K&2.22e-1&29&5.84e-1&7&1.07e0&2&2.91e-1&10&6.47e-1&5&1.12e0&4&6.92e-1&\bf3.11&\bf 2.37 \\ 
	  & &3rd&&2.5K&3.95e-1&41&1.09e1&7&2.61e0&2&5.74e-1&13&1.33e0&5&2.65e0&4&1.00e0&\bf2.53&\bf1.74 \\ 
	  & &4th&&2.5K&4.65e-1&57&1.90e0& 13& 3.90e0 & 2&8.83e-1&13&2.58e0&6&4.77e0&5&1.00e0&\bf2.15&1.13\\
	  & &4th&&2.2K&3.06e-1&35&1.24e0&7&2.46e0&2&6.06e-1&12&1.60e0&5&3.23e0&4&6.52e-1&\bf2.04&1.07 \\
	 \cline{2-20}
	 &\multirow{3}{*}{Engine}
	 &6th&\multirow{3}{*}{400K} &3.8K&5.27e0 &601&6.05e0&268&4.20e0&16&2.00e0&157&2.88e0&24&3.92e0&30&3.24e0&0.77&\bf1.62\\ 
	 &&9th  &&2.8K&1.31e0&248&4.17e-1&26&7.25e-1&3&5.43e-1&80&4.04e-1&22&6.56e-1&15&1.42e0&\bf3.41&\bf 3.51\\ 
	 &&13th&&2.5K&1.34e0&248&4.14e-1&26&7.24e-1&3&3.93e-1&54&4.91e-1&18&8.37e-1&17&9.60e-1&\bf2.31&\bf2.44\\  
	 \hline
	\multirow{7}{*}{Stiffness}&\multirow{5}{*}{Beam}
	 &1st&\multirow{5}{*}{300K}&1.9K&x&x&5.67e-1&13&9.57e-1&4&5.14e-1&63&8.95e-1&14&1.60e0&7&4.38e-1&0.77&0.85\\
	 &&2nd&&1.9K&1.31e0&358&5.45e-1&7&9.44e-1&2&4.09e-1&30&8.37e-1&10&1.36e0&4&4.50e-1&0.83&1.10\\
	 &&2nd&&1.9K&x &x &4.88e-1 & 10&9.28e-1&4&4.46e-1&60&7.84e-1&14&1.40e0&5&4.21e-1&0.86&0.94\\
	 &&3rd &&1.9K&6.67e-1&185&4.44e-1& 6&9.02e-1&2& 3.16e-1 &27&8.17e-1&10&1.44e0&4&4.03e-1&0.91&1.27\\
	 &&3rd &&1.9K& 1.19e0 &369 &4.64e-1& 8 & 9.76e-1 & 2 &3.84e-1&29&7.64e-1&11&1.48e0&4&4.57e-1&0.98&1.19\\	 
	  \cline{2-20}
	&\multirow{3}{*}{CHead}
	  &5th&2.3M&24K&3.84e2 & 829 & x & x & 1.50e2 & 22 &x&x&8.27e1&125&9.50e1&103&8.25e2&\bf5.5&\bf9.98\\ \cline{3-20}
	  &&2nd &\multirow{2}{*}{330K} &4.8K & 4.69e0 & 265 & 3.81e0 & 24 & 1.07e1 & 3 &2.45e0&112&2.97e0&92&4.54e0&94&6.56e0&\bf1.72&\bf2.68\\ 
	  &&4th &&2.6K & 4.61e-1 & 88 & 9.85e-1 & 22 & 3.17e0 & 5 &4.16e-1&26&1.03e0&18&1.74e0&17&1.06e0&\bf2.30&\bf2.54\\ 	
	  \hline
	\end{tabular}}
\caption{ Summary of solver speed for various benchmark cases. All timings are measured in seconds. The GMRES accuracy and maximum number of iterations was set to $10^{-10}$ and $1000$ respectively for all cases. The letters `x' depicts that the solver did not converge within $1000$ iterations. All LU timings are obtained using Eigen's~\cite{Eigen} partial pivoting LU solver. Level indicates the level of the dense frontal matrix in the sparse elimination tree. `T' and `I' refer to the total solve time and the number of iterations in the iterative solver respectively. Iterative solver times depicts total solve time for the iterative solver with a fast direct BDLR (ACA) solver preconditioner (low-rank computation, direct solve, iteration, etc). For BDLR, we used a depth of $1$, $3$ and $5$ for tolerances $10^{-1}$, $10^{-3}$ and $10^{-5}$ respectively. For the 4.8K and 23K cylinder head matrices, the results on the last BDLR column were obtained using a tolerance of $10^{-4}$ and a depth of $10$. We have calculated the speedups by comparing the runtime of the conventional LU solver to the lowest runtime for each case.}	
\label{table:summaryHODLRNum}	
\end{table}

\section{Conclusion and Future Work}
\label{sec:summary}

To reach our final goal of constructing a fast multifrontal solver, we need to improve the slow dense solves for the frontal matrices, which we demonstrate through various benchmarks using the HODLR solver. Using block low-rank structures like the HODLR structure  significantly reduces the memory consumption of a multifrontal sparse solver. In practice, this is currently one of the main bottlenecks on existing hardware.


One of the major challenges in constructing a fast direct sparse solver is the need for a low-rank approximation scheme that works for algebraic Schur complements, as found in multifrontal solvers. We have addressed this issue by introducing the BDLR low-rank approximation scheme which, as we have shown, is a very robust algorithm when applied to such matrices. The major advantage of BDLR over ACA is the fact that it allows us to only keep track of certain rows an columns in the multifrontal extend-add process. As we will show in a follow-up paper, this significantly improves the runtime and memory consumption of the solve process.

One of the drawbacks of this method is that it relies on off-diagonal blocks being low-rank, which may not be true asymptotically as $N \to \infty$ for points distributed on a 2D or 3D manifold. More precisely, in 2D, the rank stays fairly constant but in 3D, there is a slow growth like $n^{1/2}$ where $n$ is the size of the dense matrix or front. However, the HODLR scheme has a computational cost of $\mathcal{O}(r^2 n\log^2n)$ with a relatively small constant in front of $r^2 n\log^2n$. As a result, even if the rank $r$ increases, the method stays competitive. 

A major advantage of the mentioned direct solve algorithm is its parallel scalability. Since we make two independent recursive calls on each of the diagonal blocks, this method will scale very well. As a result, despite its somewhat higher number of flops (compared to an optimal $\mathcal{O}(n)$ method), the algorithm may run faster on large-scale parallel computers where communication and concurrency are key.

We estimate that the direct solver presented here starts being faster as soon as the rank $r$ is $r \sim 0.4 n$ compared to an $(2/3) n^3$ LU factorization algorithm. Recent work by Ambikasaran and Darve~\cite{ambikasaran2014ifmm} has overcome the growth of rank in all dimensions by requiring a compression of well-separated clusters only. It is also worth mentioning that Ho and Ying~\cite{ho2013hierarchicalA} attempt to reduce the rank, when using interpolative decompositions to build low-rank approximations, with a scheme that, in essence, is able to reduce the dimensionality of the problem in a recursive manner.

\section{Acknowledgments}
The authors would like to acknowledge Prof.\ Charbel Farhat and Dr.\ Philip Avery for providing us with the FETI test matrices. We also want to thank Prof.\ Pierre Ramet and Dr.\ Mathieu Faverge for their collaboration on this work.

Part of this research was done at Stanford University, and was supported in part by the U.S. Army Research Laboratory, through the Army High Performance Computing Research Center, Cooperative Agreement W911NF-07-0027. This material is also based upon work supported by the Department of Energy National Nuclear Security Administration under Award Number DE-NA0002373-1.

\bibliographystyle{elsarticle-harv}

\bibliography{BlackBox_HODLR,biblio_self,biblio_rest}
\end{document}